\documentclass[11pt]{article}
\usepackage{amscd}
\usepackage{amsmath}
\usepackage{amsfonts}

\newcommand{\ga}{\gamma}

\newcommand{\grad}{\mathrm{grad}}

\newcommand{\ov}[1]{\overline{#1}}

\newcommand{\ti}[1]{\tilde{#1}}
\newcommand{\tr}[2]{\textrm{tr}_{#1} #2}
\newcommand{\cur}{\mathcal{R}}
\newenvironment{proof1}{\noindent\textbf{Proof} }{\hfill$\square$\medskip}
\begin{document}
\newcounter{remark}
\newcounter{theor}
\setcounter{remark}{0}
\setcounter{theor}{1}
\newtheorem{claim}{Claim}
\newtheorem{theorem}{Theorem}[section]
\newtheorem{proposition}{Proposition}[section]
\newtheorem{conjecture}{Conjecture}[section]
\newtheorem{lemma}{Lemma}[section]
\newtheorem{defn}{Definition}[theor]
\newtheorem{corollary}{Corollary}[section]
\newenvironment{proof}[1][Proof]{\begin{trivlist}
\item[\hskip \labelsep {\bfseries #1}]}{\end{trivlist}}
\newenvironment{remark}[1][Remark]{\addtocounter{remark}{1} \begin{trivlist}
\item[\hskip
\labelsep {\bfseries #1  \thesection.\theremark}]}{\end{trivlist}}
\setlength{\arraycolsep}{2pt}
\centerline{\bf TAMING SYMPLECTIC FORMS}
\centerline{\bf AND THE CALABI-YAU EQUATION
%\footnote{Research supported in part by National Science Foundation grants DMS-0504285 and DMS-03-06600
}
\addtocounter{section}{1}
\bigskip

\centerline{{\bf Valentino Tosatti}$\,^*$, {\bf Ben Weinkove}$\,^*$ {\bf and Shing-Tung Yau}$\,^*$}

\bigskip

\centerline{$\,^*$Mathematics Department}
\centerline{Harvard University}
\centerline{Cambridge, MA 02138}

\bigskip
\bigskip

\noindent
{\bf Abstract.}  We study the Calabi-Yau equation on symplectic manifolds.  We show that Donaldson's conjecture on estimates for this equation in terms of a taming symplectic form can be reduced to an integral estimate of a scalar potential function.  Under a positive curvature condition, we show that the conjecture holds.

\bigskip
\bigskip
\noindent
{\bf 1.  Introduction}
\bigskip

Calabi's conjecture \cite{Ca1}, proved thirty years ago by the third author \cite{Ya}, states that any 
representative of the first Chern class of a compact K\"ahler manifold $(M, \omega)$ can be uniquely represented as the Ricci curvature of a K\"ahler metric in a fixed cohomology class. 
%This theorem, and its many corollaries, have laid the foundations for much of the study of K\"ahler %geometry over the last three decades.   
This can be restated in terms of volume forms as follows.
For any volume form $\sigma$  satisfying $\int_X \sigma = \int_X \omega^n$, there exists a unique K\"ahler form $\tilde{\omega}$ in $[\omega]$ solving
\begin{equation} \label{eqnCYK}
\ti{\omega}^n = \sigma,
\end{equation}
where $n$ the complex dimension of the manifold.  We call (\ref{eqnCYK}) the Calabi-Yau equation.

Recently, Donaldson \cite{D} has described how the Calabi-Yau theory could be generalized in a natural way in the setting of two-forms on four-manifolds.  His program, if carried out, would lead to many new and exciting results in symplectic geometry.  
A necessary element of this program is to obtain estimates for the Calabi-Yau equation on symplectic four-manifolds with a compatible but non-integrable almost complex structure.   The second author has recently shown in this case that the key estimates of \cite{Ya} can all be reduced to a $C^0$ estimate of a potential function, and that the equation can be solved when the Nijenhuis tensor is small in a certain sense \cite{We2}.   In this paper we will make some further progress towards Donaldson's program by showing, in a more general setting than in \cite{We2}, that the estimates for (\ref{eqnCYK}) can be reduced to an integral bound of the potential function, and that all the estimates indeed hold under a curvature assumption. 

Before stating the results precisely, we will recall some basic terminology.  An almost-K\"ahler manifold is a symplectic manifold $(M, \omega)$ together with a compatible almost complex structure $J$, meaning that $\omega$ and $J$ satisfy the two conditions
\begin{eqnarray} 
\omega(X, JX) & > & 0, \quad \textrm{for all } X \neq 0 \label{eqntame} \\
\label{eqnJinv}
\omega(JX, JY) & =&  \omega(X,Y), \quad \textrm{for all } X, Y.
\end{eqnarray}
Associated to this data is a Riemannian metric $g$ given by $g(X,Y) = \omega (X, JY)$.  We call $\omega$ an almost-K\"ahler form, and $g$ an almost-K\"ahler metric.
On the other hand, if the first condition (\ref{eqntame}) holds, but not necessarily the second (\ref{eqnJinv}), then 
we say that $\omega$ tames $J$. In this case, we can still define a Riemannian metric $g$ by
$$g(X,Y) = \frac{1}{2} \left( \omega(X, JY) + \omega(Y, JX)\right).$$
Observe that  $g$ is an almost-Hermitian metric, meaning that $g(JX,JY)=g(X,Y)$ for all vectors $X$ and $Y$.

In \cite{D}, Donaldson made the following conjecture.

\begin{conjecture} \label{conj} Let $M$ be a compact 4-manifold equipped with an almost complex structure $J$ and a taming symplectic form $\Omega$.  Let $\sigma$ be a smooth volume form on $M$ with $\int_M \sigma = \int_M \Omega^2$.  Then if $\ti{\omega}$ is an almost-K\"ahler form with $[\ti{\omega}] = [\Omega]$ and solving the Calabi-Yau equation
\begin{equation} \label{eqnCYY}
\ti{\omega}^2 = \sigma,
\end{equation}
there are $C^{\infty}$ a priori bounds on $\ti{\omega}$ depending only on $\Omega$, $J$ and $\sigma$.
\end{conjecture}

If this conjecture were to hold, it would imply, by the arguments of \cite{D} (see also the description in \cite{We2}), the following result.

\begin{conjecture} \label{conj2}  Let $M$ be a compact 4-manifold with $b^+(M)=1$ and let $J$ be an almost complex structure on $M$.  If there exists a symplectic form on $M$ taming $J$ then 
 there exists a symplectic form compatible with $J$.
\end{conjecture}

Moreover, Conjecture \ref{conj} would also imply a Calabi-Yau theorem on almost-K\"ahler 4-manifolds $(M, \omega)$ with $b^+(M)=1$:  given a normalized volume form $\sigma$ there would exist a unique almost-K\"ahler form $\ti{\omega} \in [\omega]$ satisfying $\ti{\omega}^2 = \sigma$.
For other applications of Conjecture \ref{conj}, and to see how it relates to Donaldson's broader program, see \cite{D}.  

We now state our results.  Our first result says that, in any dimension, all the \emph{a priori} bounds for Conjecture \ref{conj} can be reduced to an integral bound of a scalar potential function.  Namely, given any symplectic form $\Omega$ and almost-K\"ahler form $\ti{\omega}$ with $[\ti{\omega}] = [\Omega]$, define a smooth real-valued function $\varphi$ by 
\begin{equation} \label{eqnphi1}
\frac{1}{2n} \ti{\Delta} \varphi = 1 -  \frac{\ti{\omega}^{n-1} \wedge \Omega}{\ti{\omega}^n}, \quad \sup_M \varphi=0,
\end{equation}
where $\ti{\Delta}$ is the usual Laplacian on functions associated to the almost-K\"ahler metric $\ti{g}$.
Then we have the following result.

\bigskip
\pagebreak[3]
\noindent
{\bf Theorem 1} \ \emph{Let $\alpha>0$ be given.  Let $M$ be a compact $2n$-manifold equipped with an almost complex structure $J$ and a taming symplectic form $\Omega$.  Let $\sigma$ be a smooth volume form on $M$ with $\int_M \sigma = \int_M \Omega^n$.  Then if $\ti{\omega}$ is an almost-K\"ahler form with $[\ti{\omega}] = [\Omega]$ and solving the Calabi-Yau equation
\begin{equation} \label{eqnCYY2}
\ti{\omega}^n = \sigma,
\end{equation}
there are $C^{\infty}$ a priori bounds on $\ti{\omega}$  depending only on $\Omega$, $J$, $\sigma$, $\alpha$ and 
\begin{equation*}
I_{\alpha}(\varphi) :=\int_M e^{-\alpha \varphi} \Omega^n,
\end{equation*}
for $\varphi$ defined by (\ref{eqnphi1}).}

\bigskip

\begin{remark} \ 
\begin{enumerate}
\item The function $\varphi$ is precisely the usual K\"ahler potential in the case that $\ti{\omega}$ and $\Omega$ are K\"ahler forms, and it coincides with the `almost-K\"ahler potential' $\varphi_1$ in the terminology of \cite{We2} if they are both almost-K\"ahler.   
\item We recall a general result in K\"ahler geometry \cite{Ho}, \cite{Ti}, which is independent of the Calabi-Yau equation: the quantity $I_{\alpha}(\varphi)$
is always uniformly bounded  if $\Omega$ and $\ti{\omega}$ are K\"ahler, as long as $\alpha$ is sufficiently small (where the bounds depend only on $M$, $\Omega$, $J$ and $\alpha$).  Indeed, the supremum of all such $\alpha$ so that this quantity can be bounded independent of $\ti{\omega} \in [\Omega]$ is known as the alpha-invariant and has been much studied \cite{Ti}, \cite{TiYa}.
\item It can be easily checked that Theorem 1 still holds if $[\ti{\omega}]\neq[\Omega]$, as in \cite{We2}. Also, some of the estimates
go through if $\ti{\omega}$ is assumed to be only quasi-K\"ahler (see Section 2 for the definition).
\item  As remarked in \cite{D}, Conjecture \ref{conj2} is false in dimensions six or higher.  The deformation argument used to infer it from the first conjecture crucially uses four dimensions.  It is still possible, as far as we know, for Conjecture \ref{conj} to hold in all dimensions.  However it is also quite possible that a four dimensional argument will be needed to remove the dependence on  $I_{\alpha}(\varphi)$ in Theorem 1.
\item  Donaldson has shown, in four dimensions but in a much more general setting, that the \emph{a priori} bounds will follow if $\ti{\omega}$ is bounded in the $C^0$ norm and has fixed modulus of continuity.   In \cite{We2}, it was shown, again in four dimensions but in the case where $\Omega$ is almost-K\"ahler that the estimates can be reduced to a $C^0$ bound on $\varphi$.
\end{enumerate}
\end{remark}

Now let $g$ be the almost-Hermitian metric associated to $\Omega$ and $J$.  There exists a \emph{canonical connection} $\nabla$ associated to $(M,g, J)$.  This differs from the Levi-Civita connection, and it is described in section 2.  Under a positivity condition on the curvature of this connection, we can solve Donaldson's conjecture.  More precisely, define a tensor $$\cur_{i \ov{j} k \ov{\ell}}(g, J) = R^j_{i k \ov{\ell}} + 4 N^r_{\ov{\ell} \, \ov{j}} \ov{N^i _{\ov{r} \, \ov{k}}},$$ where $R^j_{ik \ov{\ell}}$ is the (1,1) part of the curvature of $\nabla$ and $N$ represents the Nijenhuis tensor (for precise definitions, see section 2).  We write $\cur \ge 0$ if $\cur_{i \ov{j} k \ov{\ell}} 
X^i \ov{X^j} Y^k \ov{Y^\ell} \ge 0$ for all (1,0) vectors $X$ and $Y$.

\bigskip
\noindent
{\bf Theorem 2} \ \emph{If $\cur(g, J) \ge 0$, Conjecture \ref{conj} holds.}
\bigskip

In fact under this condition we can prove Conjecture \ref{conj} in any dimension $2n$. Note that 
if $g$ were K\"ahler and the bisectional curvature of $g$ positive, then we would  have  $\cur>0$.  Hence the condition holds on $\mathbb{CP}^n$  if the pair $(\Omega, J)$ is not too far from the Fubini-Study symplectic form paired with the standard complex structure.

It will be convenient to reformulate Donaldson's conjecture as follows.  
 Let $\ti{g}$ be an almost-K\"ahler metric with K\"ahler form $\ti{\omega}$ satisfying (\ref{eqnCYK}).  Write $\sigma/n! = e^F dV_g$ where $dV_g$ is the volume form associated to $g$ and $F$ is a smooth function on $M$.  Then (\ref{eqnCYK}) can be written locally as
\begin{equation} \label{eqnCY2}
\det \ti{g} = e^{2F} \det{g},
\end{equation}
Finding bounds on $\ti{g}$ depending only on $g$, $J$ and $F$ is equivalent to solving the conjecture.

A key tool in this paper is the use of the canonical connection and the formalism of moving frames instead of the Levi-Civita connection and a normal coordinate system.   This simplifies and improves many of the estimates in \cite{We2}.  In section 2, we give the background on almost-Hermitian metrics and canonical connections and prove a number of formulas for later use.  In section 3, we prove an estimate of the metric $\ti{g}$ in terms of the potential $\varphi$.  This is the analogue of the second order estimate of \cite{Ya}.  In section 4, we give an estimate of the first derivative of $\ti{g}$ in terms of $\ti{g}$ itself.  This is the analogue of the well-known third order estimate of \cite{Ya} (which was inspired by that of Calabi \cite{Ca2}).  In fact, this section is not strictly necessary to complete the proofs of the theorems, since we could have instead appealed to the argument in \cite{We2}, adapting the technique of Evans \cite{Ev} and Krylov \cite{Kr}.  However, we have included the `third order' estimate since it is self-contained and in the spirit of the rest of the paper.  In section 5, we give a proof of Theorem 1.  We make use of a Moser iteration argument as in \cite{Ya}, but applied to the exponential of $\varphi$ in a similar way to \cite{We1}, to obtain a $C^0$ estimate of $\varphi$ depending on $I_{\alpha}(\varphi)$. We also prove higher order estimates using a bootstrapping argument.  Finally, in section 6, we give a proof of Theorem 2.

\setcounter{equation}{0}
\addtocounter{section}{1}
\setcounter{remark}{0}
\bigskip
\bigskip
\noindent
{\bf 2. Almost-Hermitian manifolds and the canonical connection}
\bigskip

In this section, we give some background on almost-Hermitian manifolds, almost- and quasi-K\"ahler manifolds,  the canonical connection and its torsion and curvature. 
Many of the results of this section are well-known, and so we have omitted the proofs of several of the lemmas (a good reference for this material is \cite{Ga}).  On the other hand, whenever  precise formulas do not seem to be readily available in the literature, we have provided the arguments.

\bigskip \pagebreak[3]
\noindent
{\bf 2.1. Almost-Hermitian metrics and connections}
\bigskip

Let $M$ be a manifold of dimension $2n$ with an almost complex structure $J$ and a Riemannian metric $g$ satisfying
$$g(J X, JY) = g(X,Y),$$
for all tangent vectors $X$ and $Y$.  We say that $(M,J,g)$ is an \emph{almost-Hermitian manifold}.

Write $T^{\mathbb{R}}_p M$ for the (real) tangent space of $M$ at a point $p$.  In the following we will drop the subscript $p$.  Denote the complexified tangent space by $T^{\mathbb{C}}M = T^{\mathbb{R}} M \otimes \mathbb{C}$.  Extending $g$ and $J$ linearly to $T^{\mathbb{C}}M$, we see that the complexified tangent space can be decomposed as
$$T^{\mathbb{C}} M = T' M \oplus T''M,$$
where $T'M$ and $T''M$ are the eigenspaces of $J$ corresponding to eigenvalues $\sqrt{-1}$ and $- \sqrt{-1}$ respectively.  Extending $J$ to forms, we can uniquely decompose $m$-forms into $(p,q)$-forms for each $p$, $q$ with $p+q=m$.

Choose a local unitary frame $\{ e_1, \ldots, e_n \}$ for $T'M$ with respect to the Hermitian inner product induced from $g$, and let $\{ \theta^1, \ldots, \theta^n \}$ be a dual coframe.  The metric $g$ can be written as
$$g =   \theta^i \otimes \ov{\theta^i} + \ov{\theta^i} \otimes \theta^i,$$
where here, and henceforth, we are summing over repeated indices.

Let $\nabla$ be an affine connection on $T^{\mathbb{R}} M$, which we extend linearly to $T^{\mathbb{C}} M$.  We say that $\nabla$ is an \emph{almost-Hermitian connection} if
$$\nabla J = \nabla g = 0.$$
It is well-known that (see e.g. \cite{KN}):

\begin{lemma} Almost-Hermitian connections always exist on almost-Hermitian manifolds.
\end{lemma}

From now on, assume that $\nabla$ is almost-Hermitian.  Observe that for $i=1, \ldots, n$,
$$J ( \nabla e_i) = \sqrt{-1} \nabla e_i,$$ 
and hence $\nabla e_i \in T'M \otimes (T^{\mathbb{C}}(M))^*$.  Then locally there exists a matrix of complex valued 1-forms $\{ \theta_i^j \}$, called the \emph{connection 1-forms}, such that
$$\nabla e_i = \theta_i^j e_j.$$
Applying $\nabla$ to $g(e_i, \ov{e_j})$ and using the condition $\nabla g=0$ we see that $\{ \theta_i^j \}$ satisfies the skew-Hermitian property
$$\theta_i^j + \ov{\theta_j^i} = 0.$$
Now define the \emph{torsion} $\Theta = (\Theta^1, \ldots, \Theta^n)$ of $\nabla$ by
\begin{equation} \label{eqnstructure1}
d \theta^i = - \theta_j^i \wedge \theta^j + \Theta^i, \qquad \textrm{for } i=1, \ldots, n.
\end{equation}
Notice that the $\Theta^i$ are 2-forms.  Equation (\ref{eqnstructure1}) is known as the \emph{first structure equation}.
Define the \emph{curvature} $\Psi = \{ \Psi_j^i \}$   of $\nabla$ by
\begin{equation} \label{eqnstructure2}
d \theta_j^i = - \theta_k^i \wedge \theta_j^k + \Psi_j^i.
\end{equation}
Note that $\{ \Psi_j^i \}$ is a  skew-Hermitian matrix of 2-forms.  Equation (\ref{eqnstructure2}) is known as the \emph{second structure equation}.

\bigskip \pagebreak[3]
\noindent
{\bf 2.2. The canonical connection}
\bigskip

We have the following lemma (see e.g. \cite{Ga}).

\begin{lemma}
There exists a unique almost-Hermitian connection $\nabla$ on $(M,J,g)$ whose torsion $\Theta$ has everywhere vanishing $(1,1)$ part.
\end{lemma}

Such a connection is known as the \emph{second canonical connection} and was first introduced by Ehresmann and Libermann in \cite{EL}.  It is also sometimes referred to as the \emph{Chern connection}, since when $J$ is integrable it coincides with the connection defined in \cite{Ch}.
We will call it simply the \emph{canonical connection}.  

Define functions $T^i_{jk}$ and  $N^i_{\ov{j} \, \ov{k}}$ by
\begin{eqnarray*}
(\Theta^i)^{(2,0)} & = & T^i_{jk} \theta^j \wedge \theta^k \\
(\Theta^i)^{(0,2)} & = & N^i_{\ov{j} \, \ov{k}} \ov{\theta^j} \wedge \ov{\theta^k},
\end{eqnarray*}
with $T^i_{jk} = - T^i_{kj}$ and $N^i_{\ov{j} \, \ov{k}} = - N^i_{\ov{k} \, \ov{j}}$.

\begin{lemma} \label{lemmanijenhuis}  The (0,2) part of the torsion is independent of the choice of metric.
\end{lemma}

Indeed $(\Theta^i)^{(0,2)}$ can be regarded as the Nijenhuis tensor of $J$. For a proof of this lemma, see section 3.

Let's consider now the real $(1,1)$ form 
$$\omega=\sqrt{-1} \theta^i\wedge\ov{\theta^i}.$$
We say that $(M,J,g)$ is \emph{almost-K\"ahler} if $d\omega=0$, and that it is \emph{quasi-K\"ahler} if
$(d\omega)^{(1,2)}=0$. An almost-K\"ahler or quasi-K\"ahler manifold with $J$ integrable is a K\"ahler manifold.  Observe from the first structure equation,
\begin{eqnarray*}
d\omega&=& \sqrt{-1}(\Theta^i\wedge\ov{\theta^i}-\theta^i\wedge\ov{\Theta^i}) \\
& = &
\sqrt{-1}(N^i_{\ov{j}\,\ov{k}}\ov{\theta^i}\wedge\ov{\theta^j}\wedge\ov{\theta^k}-\ov{N^i_{\ov{j}\,\ov{k}}}\theta^i\wedge\theta^j\wedge\theta^k
\\&& \mbox{} +T^i_{jk}\ov{\theta^i}\wedge\theta^j\wedge\theta^k-\ov{T^i_{jk}}\theta^i\wedge\ov{\theta^j}\wedge\ov{\theta^k}).
\end{eqnarray*}
Thus we have the following alternative definitions using the torsion of the canonical connection.

\begin{lemma} An almost-Hermitian manifold $(M,J,g)$ is almost-K\"ahler if and only if
$$T^i_{jk}=0,$$
\begin{equation}\label{bianchin}
N_{\ov{i}\,\ov{j}\,\ov{k}}+N_{\ov{j}\,\ov{k}\,\ov{i}}+N_{\ov{k}\,\ov{i}\,\ov{j}}=0,
\end{equation}
where $N_{\ov{i}\,\ov{j}\,\ov{k}}=N^i_{\ov{j}\,\ov{k}}$, and is quasi-K\"ahler if and only if
$$T^i_{jk}=0.$$
\end{lemma}

In particular on a quasi-K\"ahler manifold the torsion of the canonical connection has only a
$(0,2)$ component
$$\Theta^i=N^i_{\ov{j}\,\ov{k}}\ov{\theta^j}\wedge\ov{\theta^k}.$$

\bigskip \pagebreak[3]
\noindent
{\bf 2.3. Curvature identities}
\bigskip

Let $(M,J,g)$ be an almost-Hermitian manifold and let $\nabla$ be the canonical connection with torsion $\Theta$ and curvature  $\Psi$.  Define $R_{ik \ov{\ell}}^j$,  $K^i_{jk\ell}$ and $K^i_{j \ov{k} \, \ov{\ell}}$ by
\begin{eqnarray*}
(\Psi_i^j)^{(1,1)} & = & R_{ik \ov{\ell}}^j \theta^k \wedge \ov{\theta^\ell} \\
(\Psi^i_j)^{(2,0)} & = & K^i_{jk\ell} \theta^k \wedge \theta^\ell\\
(\Psi^i_j)^{(0,2)} & = & K^i_{j \ov{k} \, \ov{\ell}} \ov{\theta^k} \wedge \ov{\theta^\ell},
\end{eqnarray*}
with $K^i_{jk\ell} = - K^i_{j\ell k}$ and  $K^i_{j\ov{k} \, \ov{\ell}} = - K^i_{j\ov{\ell} \, \ov{k}}$.  We define the \emph{Ricci curvature} and \emph{scalar curvature} of the canonical connection to be the tensors $R_{k\ov{\ell}} = R^i_{ik\ov{\ell}}$ and $R = R_{k \ov{k}}$ respectively.

We will now derive some curvature identities.  Applying the exterior derivative to the first and second structure equations, we obtain the \emph{first Bianchi identity},
\begin{equation} \label{eqnfirstBianchi}
d \Theta^i = \Psi^i_j \wedge \theta^j - \theta^i_j \wedge \Theta^j,
\end{equation}
and \emph{second Bianchi identity},
\begin{equation} \label{eqnsecondBianchi}
d \Psi^i_j = \Psi^i_k \wedge \theta^k_j - \theta^i_k \wedge \Psi^k_j.
\end{equation}

Let us rewrite these.  First, define
$T^i_{jk,p}$, $T^i_{jk,\ov{p}}$ by
\begin{equation} \label{derivativeT}
dT^i_{jk} + \theta^i_p T^p_{jk} - T^i_{pk} \theta^p_j - T^i_{jp} \theta^p_k 
 =  T^i_{jk, p} \theta^p + T^i_{jk, \ov{p}} \ov{\theta^p},
\end{equation}
and $N^i_{\ov{j} \, \ov{k}, p}$ and $N^i_{\ov{j} \, \ov{k}, \ov{p}} $ by
\begin{equation}\label{derivativeN}
dN^i_{\ov{j} \, \ov{k}} + \theta^i_p N^p_{\ov{j} \, \ov{k}} - N^i_{\ov{p} \, \ov{k}} \ov{\theta^p_j} - N^i_{\ov{j} \ov{p}} \ov{\theta^p_k}  =  N^i_{\ov{j} \, \ov{k}, p} \theta^p + N^i_{\ov{j} \, \ov{k}, \ov{p}} \ov{\theta^p}.
\end{equation}
Then the first Bianchi identity can be written as
\begin{eqnarray*}
dT^i_{jk} \wedge \theta^j \wedge \theta^k - T^i_{jk} \theta^j_p \wedge \theta^p \wedge \theta^k + T^i_{jk} \Theta^j \wedge \theta^k 
+ T^i_{jk} \theta^j \wedge \theta^k_p \wedge \theta^p \\ - T^i_{jk} \theta^j \wedge \Theta^k + dN^i_{\ov{j} \, \ov{k}} \wedge \ov{\theta^j} \wedge  \ov{\theta^k} - N^i_{\ov{j} \, \ov{k}} \ov{\theta^j_p} \wedge \ov{\theta^p} \wedge \ov{\theta^k}  \\
+ N^i_{\ov{j} \, \ov{k}} \ov{\theta^j} \wedge \ov{\theta^k_p} \wedge \ov{\theta^p} + N^i_{\ov{j} \, \ov{k}} \ov{\Theta^j} \wedge \ov{\theta^k} - N^i_{\ov{j} \, \ov{k}} \ov{\theta^j} \wedge \ov{\Theta^k} \\
= K^i_{jk\ell} \theta^k \wedge \theta^\ell \wedge \theta^j + R^i_{jk \ov{\ell}} \theta^k \wedge \ov{\theta^\ell} \wedge \theta^j + K^i_{j \ov{k} \, \ov{\ell}} \ov{\theta^k} \wedge \ov{\theta^\ell} \wedge \theta^j\\
-T^j_{k\ell}\theta^i_j\wedge\theta^k\wedge\theta^\ell-N^j_{\ov{k}\,\ov{\ell}}\theta^i_j\wedge\ov{\theta^k}\wedge\ov{\theta^\ell}.
\end{eqnarray*}
After substituting from (\ref{derivativeT}) and (\ref{derivativeN}), and comparing bidegrees, we arrive at the following four identities:
\begin{eqnarray*}
(T^i_{jk,\ell}+ 2 T^i_{pj} T^p_{k\ell}- K^i_{jk\ell} ) \theta^j \wedge \theta^k \wedge \theta^\ell  & = & 0\\
(T^i_{jk, \ov{\ell}} + 2 \ov{N^p_{\ov{j} \, \ov{k}}} N^i_{\ov{p} \, \ov{\ell}} -R^i_{jk\ov{\ell}}) \theta^j \wedge \theta^k \wedge \ov{\theta^\ell}  & = & 0 \\
(2T^i_{pj} N^p_{\ov{k} \, \ov{\ell}} +N^i_{\ov{k} \, \ov{\ell}, j} - K^i_{j \ov{k} \, \ov{\ell}}) \theta^j \wedge \ov{\theta^k} \wedge \ov{\theta^\ell} & = & 0 \\
(N^i_{\ov{j} \, \ov{k} , \ov{\ell}} + 2 N^i_{\ov{p} \, \ov{j}} \ov{T^p_{k\ell}} ) \ov{\theta^j} \wedge \ov{\theta^k} \wedge \ov{\theta^\ell} & = & 0,
 \end{eqnarray*}
which are equivalent to:
\begin{eqnarray} 
\setlength{\arraycolsep}{1pt}
T^i_{jk,\ell} + T^i_{k\ell,j} + T^i_{\ell j,k}
 + 2T^i_{pj} T^p_{k\ell}+ 2T^i_{pk} T^p_{\ell j} + 2T^i_{p\ell} T^p_{jk}  
& = & K^i_{jk\ell} + K^i_{k\ell j} + K^i_{\ell jk} \qquad \quad  \\
2T^i_{jk, \ov{\ell}} + 4 \ov{N^p_{\ov{j} \, \ov{k}}} N^i_{\ov{p} \, \ov{\ell}} & = & R^i_{jk \ov{\ell}} - R^i_{kj \ov{\ell}} \label{curvatureidentities}\\
2T^i_{pj} N^p_{\ov{k} \, \ov{\ell}} + N^i_{\ov{k} \,  \ov{\ell},j} & = & K^i_{j \ov{k} \, \ov{\ell}} \label{curvatureidentities3} \\
N^i_{\ov{j} \, \ov{k}, \ov{\ell}} + N^i_{\ov{k} \, \ov{\ell}, \ov{j}}+N^i_{\ov{\ell} \, \ov{j}, \ov{k}} + 2 N^i_{\ov{p} \, \ov{j}} \ov{T^p_{k\ell}} + 2 N^i_{\ov{p} \, \ov{k}} \ov{T^p_{\ell j}} +2 N^i_{\ov{p} \, \ov{\ell}} \ov{T^p_{jk}} & = & 0.
\end{eqnarray}

By a similar reasoning, we obtain the following from the second Bianchi identity:
\begin{eqnarray*}
(K^i_{jk\ell,p} + 2T^q_{k\ell} K^i_{jqp} - R^i_{jk\ov{q}} \ov{N^q_{\ov{\ell} \, \ov{p}}} ) \theta^k \wedge \theta^\ell \wedge \theta^p & =& 0 \\
(K^i_{jk\ell, \ov{p}} - R^i_{jk \ov{p},\ell} + R^i_{jq\ov{p}} T^q_{k\ell} + 2K^i_{j\ov{q} \, \ov{p}} \ov{N^q_{\ov{k} \, \ov{\ell}}}) \theta^k \wedge \theta^\ell \wedge \ov{\theta^p}& = & 0\\
(R^i_{jk \ov{\ell}, \ov{p}} + K^i_{j \ov{\ell} \, \ov{p},k} + 2 K^i_{j q k} N^q_{\ov{\ell} \, \ov{p}} - R^i_{jk\ov{q}} \ov{T^q_{\ell p}} )\theta^k \wedge \ov{\theta^\ell} \wedge \ov{\theta^p} & = & 0 \\
(K^i_{j \ov{k} \, \ov{\ell}, \ov{p}}  + R^i_{jq\ov{k}} N^q_{\ov{\ell} \, \ov{p}} + 2K^i_{j \ov{q} \, \ov{p}} \ov{T^q_{k\ell}})       \ov{\theta^k} \wedge \ov{\theta^\ell} \wedge \ov{\theta^p} & = & 0,
\end{eqnarray*}
where $K^i_{jk\ell,p}$, $K^i_{jk\ell, \ov{p}}$ etc. are defined in the obvious way.  The above four identities can be rewritten as
\begin{eqnarray} \nonumber
\lefteqn{K^i_{jk\ell,p} + K^i_{j\ell p,k} + K^i_{jpk,\ell} + 2T^q_{k\ell} K^i_{jqp} + 2T^q_{\ell p} K^i_{jqk}}  \\
 && \mbox{} + 2T^q_{pk} K^i_{jq\ell} - R^i_{jk\ov{q}} \ov{N^q_{\ov{\ell} \, \ov{p}}} - R^i_{j\ell\ov{q}} \ov{N^q_{\ov{p} \, \ov{k}}} - R^i_{jp\ov{q}} \ov{N^q_{\ov{k} \, \ov{\ell}}}  =  0 \\
&&  2K^i_{jk\ell,\ov{p}} - R^i_{jk\ov{p},\ell} + R^i_{j\ell\ov{p},k} + 2R^i_{jq\ov{p}} T^q_{k\ell} + 4K^i_{j\ov{q} \, \ov{p}} \ov{N^q_{\ov{k} \, \ov{\ell}}}  =  0 \label{curvatureidentities2}\\
&& R^i_{jk\ov{\ell}, \ov{p}} - R^i_{jk\ov{p},\ov{\ell}} + 2K^i_{j\ov{\ell} \, \ov{p},k} + 4K^i_{jqk} N^q_{\ov{\ell}\, \ov{p}} - 2R^i_{jk\ov{q}} \ov{T^q_{\ell p}} =0 \\ \nonumber
\lefteqn{K^i_{j \ov{k} \, \ov{\ell},\ov{p}} + K^i_{j \ov{\ell} \, \ov{p},\ov{k}} + K^i_{j \ov{p} \, \ov{k} ,\ov{\ell}} + R^i_{j q \ov{k}} N^q_{\ov{\ell} \, \ov{p}} + R^i_{j q \ov{\ell}} N^q_{\ov{p} \, \ov{k}} + R^i_{jq\ov{p}}N^q_{\ov{k} \, \ov{\ell}} }\\
&& \mbox{} + 2K^i_{j\ov{q} \, \ov{p}} \ov{T^q_{k\ell}} + 2K^i_{j\ov{q} \, \ov{k}} \ov{T^q_{\ell p}} + 2K^i_{j\ov{q} \, \ov{\ell}} \ov{T^q_{pk}} =0.
\end{eqnarray}

Now assume that $(M, g, J)$ is quasi-K\"ahler, so that the (2,0) part of the torsion vanishes.  Then
 \eqref{curvatureidentities}, \eqref{curvatureidentities3} and \eqref{curvatureidentities2} above simplify to
\begin{equation}
\label{commutation1}
4 \ov{N^p_{\ov{j} \, \ov{k}}} N^i_{\ov{p} \, \ov{\ell}}  = R^i_{jk \ov{\ell}} - R^i_{kj \ov{\ell}},
\end{equation}
\begin{equation}
\label{commutation3}
N^i_{\ov{k} \,  \ov{\ell},j}  =  K^i_{j \ov{k} \, \ov{\ell}},
\end{equation}
\begin{equation}
\label{commutation2}
2K^i_{jk\ell,\ov{p}}  + 4K^i_{j\ov{q} \, \ov{p}} \ov{N^q_{\ov{k} \, \ov{\ell}}}  =  R^i_{jk\ov{p},\ell} - R^i_{j\ell\ov{p},k}.
\end{equation}
Recall that the curvature matrix $(\Psi^i_j)$ is skew-Hermitian, hence
\begin{equation}\label{curvident}
K^i_{jk\ell}=\ov{K^j_{i\ov{\ell}\,\ov{k}}},\quad R^i_{jk\ov{\ell}}=\ov{R^j_{i\ell\ov{k}}}.
\end{equation}
From this we compute
\begin{eqnarray}\label{switchr}\nonumber
R^i_{jk\ov{\ell}}&=&R^i_{kj\ov{\ell}}+4 N^i_{\ov{p}\,\ov{\ell}}\ov{N^p_{\ov{j}\,\ov{k}}}\\ \nonumber
&=&\ov{R^k_{i\ell\ov{j}}}+4 N^i_{\ov{p}\,\ov{\ell}}\ov{N^p_{\ov{j}\,\ov{k}}}\\ \nonumber
&=&\ov{R^k_{\ell i\ov{j}}}+4 N^i_{\ov{p}\,\ov{\ell}}\ov{N^p_{\ov{j}\,\ov{k}}}+4N^p_{\ov{i}\,\ov{\ell}}\ov{N^k_{\ov{p}\,\ov{j}}}\\ 
&=&R^\ell_{kj\ov{i}}+4 N^i_{\ov{p}\,\ov{\ell}}\ov{N^p_{\ov{j}\,\ov{k}}}+4N^p_{\ov{i}\,\ov{\ell}}\ov{N^k_{\ov{p}\,\ov{j}}},
\end{eqnarray}
giving us the following formula for the Ricci curvature
\begin{equation}\label{riccicurv}
R_{k\ov{\ell}}=R^i_{ik\ov{\ell}}= R^\ell_{ki\ov{i}}+4 N^i_{\ov{p}\,\ov{\ell}}\ov{N^p_{\ov{i}\,\ov{k}}}+4N^p_{\ov{i}\,\ov{\ell}}\ov{N^k_{\ov{p}\,\ov{i}}}.
\end{equation}

%\begin{equation}\label{derivativeN}
%dN^i_{\ov{j} \, \ov{k}} + \theta^i_p N^p_{\ov{j} \, \ov{k}} - N^i_{\ov{p} \, \ov{k}} \ov{\theta^p_j} - N^i_{\ov{j} \ov{p}} %\ov{\theta^p_k}  =  N^i_{\ov{j} \, \ov{k}, p} \theta^p + N^i_{\ov{j} \, \ov{k}, \ov{p}} \ov{\theta^p} 
%\end{equation}
%\begin{equation*}
%dN^i_{\ov{j} \, \ov{k}, p} + \theta^i_q N^q_{\ov{j} \ov{k}, p} - N^i_{\ov{q} \, \ov{k}, p} \ov{\theta^q_j} - N^i_{\ov{j} \, \ov{q}, p} \ov{\theta^q_k} - N^i_{\ov{j} \, \ov{k}, q} \theta^q_p  =  N^i_{\ov{j} \, \ov{k},pq}\theta^q + N^i_{\ov{j} \, \ov{k}, p\ov{q}} \ov{\theta^q}.
%\end{equation*}
%On an almost-Hermitian manifold we have the following identities:
%\begin{align} 
%& K^i_{jkl} + K^i_{klj} + K^i_{ljk} =0 \label{eqnbianchi1} \\
%& K^i_{jkl,\ov{p}} + K^i_{klj, \ov{p}} + K^i_{ljk, \ov{p}} =0 \label{eqnbianchi2} \\
%& 2 K^i_{jkl, \ov{p}} - R^i_{jk\ov{p},l} + R^i_{jl \ov{p},k} =0 \label{eqncurv1} \\
%& N^i_{\ov{k} \, \ov{l}, j} =  \ov{K^j_{ilk}} \label{eqncurv2} \\
%& N^i_{\ov{k} \, \ov{l}, jp} = \ov{K^j_{ilk, \ov{p}}}. \label{eqncurv3}
%\end{align}
%From the above we also have the following commutation formula
%\begin{eqnarray} \label{eqncomm1}
%N^i_{\ov{j} \, \ov{k}, rs} - N^i_{\ov{j} \, \ov{k}, sr} & = & 2 ( N^i_{\ov{j} \, \ov{k}, \ov{p}} \ov{N^p_{\ov{r} \, \ov%{s}}} - K^i_{prs} N^p_{\ov{j} \, \ov{k}} + N^i_{\ov{p}\, \ov{k}} \ov{K^p_{j \ov{r} \, \ov{s}}} + N^i_{\ov{j} \, \ov{p}} \ov{K^p_{k \ov{r} \, \ov{s}}}) \quad 
%\end{eqnarray}

\bigskip \pagebreak[3]
\noindent
{\bf 2.4. The canonical Laplacian}
\bigskip

Suppose that $(M,J,g)$ is almost-Hermitian and let $\nabla$ be its canonical connection. Let $f$ be a function on $M$.  We define the \emph{canonical Laplacian} $\Delta$ of $f$ by
$$\Delta f  = \sum_i \left( (\nabla \nabla f) (e_i, \ov{e_i}) + (\nabla \nabla f) (\ov{e_i}, e_i) \right).$$
This expression is independent of the choice of unitary frame.  

Define $f_i$ and $f_{\ov{i}}$ by
\begin{equation} \label{eqndu}
df = f_i \theta^i + f_{\ov{i}} \theta^i.
\end{equation}
Writing $\partial f$ and $\ov{\partial} f$ for the (1,0) and (0,1) parts of $df$ respectively we see that $\partial f = f_i \theta^i$ and $\ov{\partial} f = f_{\ov{i}} \theta^i$.
Applying the exterior derivative to (\ref{eqndu}) and using the first structure equation we obtain
\begin{eqnarray} \nonumber
0 & = & d f_i \wedge \theta^i - f_i \theta_j^i \wedge \theta^j + f_i \Theta^i + d f_{\ov{i}} \wedge \theta^i  - f_{\ov{i}} \, \ov{\theta_j^i} \wedge \ov{\theta^j} + f_{\ov{i}} \ov{\Theta^i} \\ \label{eqnddu}
& = & (df_i - f_j \theta_i^j) \wedge \theta^i + (df_{\ov{i}} - f_{\ov{j}} \ov{\theta_i^j} ) \wedge \ov{\theta^i} + f_i \Theta^i + f_{\ov{i}} \ov{\Theta^i}.
\end{eqnarray}
Define $f_{ik}$, $f_{i\ov{k}}$, $f_{\ov{i} k}$ and $f_{\ov{i} \, \ov{k}}$ by
\begin{align*}
& df_i - f_j \theta_i^j = f_{ik} \theta^k + f_{i \ov{k}} \ov{\theta^k} \\
& df_{\ov{i}} - f_{\ov{j}} \ov{\theta_i^j} = f_{\ov{i} k} \theta^k + f_{\ov{i} \, \ov{k}} \ov{\theta^k}.
\end{align*}
Taking the (1,1) part of (\ref{eqnddu}) we see that
$$f_{i \ov{k}} \ov{\theta^k} \wedge \theta^i + f_{\ov{i} k} \theta^k \wedge \ov{\theta^i} = 0,$$
and hence
$$f_{i \ov{k}} = f_{\ov{k} i}.$$

Now calculate
\begin{eqnarray*}
\nabla \nabla f & = & \nabla (f_i \theta^i + f_{\ov{i}} \ov{\theta^i}) \\
& = & df_i \otimes \theta^i - f_i \theta_j^i \otimes \theta^j + d f_{\ov{i}} \otimes \ov{\theta^i} - f_{\ov{i}} \ov{\theta_j^i} \otimes \ov{\theta^j} \\
& = & (f_{i j} \theta^j + f_{i \ov{j}} \ov{\theta^j} ) \otimes \theta^i + (f_{\ov{i} j} \theta^j + f_{\ov{i} \, \ov{j}} \ov{\theta^j}) \otimes \ov{\theta^i}.
\end{eqnarray*}
Hence
\begin{equation} \label{eqnDeltau}
\Delta f =  f_{i \ov{i}} + f_{\ov{i} i} = 2  f_{i \ov{i}}.
\end{equation}
There are other ways of writing $\Delta f$.  
\begin{lemma} \label{lemmalap}
\begin{eqnarray} \label{eqnlemma1}
\Delta f  & = & - 2 \sum_{i}  (d \partial f)^{(1,1)}(e_i, \ov{e_i}) \\ \label{eqnlemma2}
& = & 2 \sum_{i}  (d \ov{\partial} f)^{(1,1)}(e_i, \ov{e_i}) \\ \label{eqnlemma3}
& = &  \sqrt{-1} \sum_{i}  (d (J d f))^{(1,1)}(e_i, \ov{e_i}),
\end{eqnarray}
where $J$ acts on a 1-form $\alpha$ by $(J\alpha)(X) = \alpha (J (X))$ for a vector $X$.
\end{lemma}
\begin{proof1}
Calculate
\begin{eqnarray}\label{derivf}\nonumber
d \partial f & = & d(f_i \theta^i) \\ \nonumber
& = & (f_{ik} \theta^k + f_{i \ov{k}} \ov{\theta^k} + f_j \theta_i^j) \wedge \theta^i - f_i \theta_j^i \wedge\theta^j + f_i \Theta^i\\
& = & f_{ik} \theta^k \wedge \theta^i+ f_{i \ov{k}} \ov{\theta^k}  \wedge \theta^i  + f_i \Theta^i.
\end{eqnarray}
Hence 
\begin{eqnarray} \label{eqndpartialf}
(d \partial f)^{(1,1)} = - f_{i \ov{k}} \theta^i \wedge \ov{\theta^k},
\end{eqnarray}
and (\ref{eqnlemma1}) follows from (\ref{eqnDeltau}).  For (\ref{eqnlemma2}), just observe that $\partial = d - \ov{\partial}$ and $d^2=0$.  For (\ref{eqnlemma3}), recall that $J \theta^i = \sqrt{-1} \theta^i$.  Then
\begin{eqnarray*}
d(J df) & = & d( J (f_i \theta^i + f_{\ov{i}} \ov{\theta^i} ))\\ 
& = & \sqrt{-1} d (f_i \theta^i - f_{\ov{i}} \ov{\theta^i}) \\
& = & \sqrt{-1} d (\partial f - \ov{\partial} f) \\
& = & 2 \sqrt{-1} d \partial f.
\end{eqnarray*}
\end{proof1}

Finally we have the following lemma.

\begin{lemma}  If the metric $g$ is quasi-K\"ahler then
 the canonical Laplacian is equal to the usual Laplacian of
the Levi-Civita connection of $g$.
\end{lemma}
\begin{proof1}
In fact, the Laplacian of the Levi-Civita connection applied to a function $f$ is given by the trace
of the map $F:TM\to TM$ defined by
$$F(X)=\nabla_X (\grad_g f)+\tau(\grad_g f,X),$$
where $\nabla$ is the canonical connection and $\tau$ is its torsion (see for example \cite{KN} p.282). But if $g$ is quasi-K\"ahler
$\tau$ is just the Nijenhuis tensor, which maps $T''M\otimes T''M\to T'M$ and so
the second term above has trace zero.

\end{proof1}

\addtocounter{section}{1}
\setcounter{lemma}{0}
\setcounter{equation}{0}
\bigskip
\bigskip \pagebreak[3]
\noindent
{\bf 3. Estimate of the metric}
\bigskip

In this section we will prove an estimate on an almost-K\"ahler metric $\ti{g}$ solving (\ref{eqnCY2}), in terms of the potential function $\varphi$.  Recall that $\varphi$ is defined by (\ref{eqnphi1}), which can be rewritten as
\begin{equation} \label{eqntiDeltaphi}
\ti{\Delta} \varphi =  2n - \tr{\ti{g}}{g},
\end{equation}
since
\begin{equation}\label{cohomolo}
\tr{\ti{g}}{g} = 2n \frac{\ti{\omega}^{n-1} \wedge \Omega}{\ti{\omega}^n}.
\end{equation}
To see (\ref{cohomolo}),  observe that
$$g_{ij} = \frac{1}{2} \left( \Omega_{ik} J_j^{\ k} + \Omega_{jk} J_i^{\ k} \right),$$
and so we have
$$\tr{\ti{g}}{g} = \ti{g}^{ij} g_{ij} = \ti{J}^{ik} \Omega_{ik},$$
where $\ti{J}^{ik} = \ti{g}^{il} J_l^{\ k}$.
Working in a coordinate system in which $\ti{\omega} = dx^1 \wedge dx^2 +\dots+ dx^{2n-1} \wedge dx^{2n}$ and $\ti{g}_{ij} = \delta_{ij}$ at a fixed point $p$ in $M$ we see that
$$\ti{J}^{ik} \Omega_{ik} = 2n \frac{\ti{\omega}^{n-1} \wedge \Omega}{\ti{\omega}^n},$$
as required.  

The estimate we wish to prove in this section is:

\begin{theorem} \label{theoremC2}
Let $\ti{g}$ be an almost-K\"ahler metric solving the Calabi-Yau equation (\ref{eqnCY2}), where $g$ is an almost-Hermitian metric.  Then
there exist constants $C$ and $A$ depending only on $J$, $R$, the lower bound of $\cur_{i \ov{j} k \ov{l}}$, $\sup|F|$ and the lower bound of $\Delta F$ such that
$$\emph{tr}_g \ti{g} \le C e^{A(\varphi - \inf_M \varphi)}.$$ 
\end{theorem}

We first compute some general formulas which are completely independent of the Calabi-Yau equation.  Let $(M, J)$ be an almost complex manifold with two almost-Hermitian metrics $g$ and $\ti{g}$.  Let $\theta^i$ and $\ti{\theta}^i$ be local unitary coframes for $g$ and $\ti{g}$ respectively.  Denote by $\nabla$ and $\ti{\nabla}$  the associated canonical connections.  We will use $\ti{\Theta}$, $\ti{\Psi}$ etc. to denote the torsion, curvature and so on with respect to $\ti{\nabla}$. 
%We may pick the frames so that at a fixed point $p$ on $M$, the connection 1-forms $\theta_i^j$ and $\ti{\theta}_i^j$ vanish (although their derivatives will not vanish, in general).  
Define local matrices $( a_j^i )$ and $(b_j^i)$ by
\begin{eqnarray} \label{eqnaji}
\ti{\theta}^i &=& a_j^i \theta^j \\ \label{eqnb}
 \theta^j & = & b_i^j \ti{\theta}^i,
 \end{eqnarray}
so that $a_j^i b_i^k = \delta_j^k$.  Define a function $u$ by 
$$u =  a_j^i \ov{a_j^i} = \frac{1}{2} \textrm{tr}_g \ti{g}.$$
Differentiating (\ref{eqnaji}) and using the first structure equations we obtain
\begin{eqnarray*}
- \ti{\theta}_k^i \wedge \ti{\theta}^k + \ti{\Theta}^i & = & da^i_j \wedge \theta^j - a_j^i \theta_k^j \wedge \theta^k + a_j^i \Theta^j.
\end{eqnarray*}
Using (\ref{eqnb}) and rearranging, we have
\begin{eqnarray} \label{eqnbda}
(b_k^j d a_j^i - a_j^i b_k^\ell \theta_\ell^j + \ti{\theta}_k^i ) \wedge \ti{\theta}^k & = & \ti{\Theta}^i - a_j^i \Theta^j.
\end{eqnarray}
Taking the $(0,2)$ part of this equation, we see that
\begin{equation}\label{eqnN0}
\ti{N}^i_{\ov{j} \, \ov{k}} =  \ov{b^r_j b^s_k} a^i_t N^t_{\ov{r} \, \ov{s}},
\end{equation}
which shows that the $(0,2)$ part of the torsion is independent of the choice of the metric (thus giving the proof of Lemma \ref{lemmanijenhuis}).

By the definition of the canonical connection, the right hand side of \eqref{eqnbda} has no (1,1)-part.  Hence there exist functions $a_{k\ell}^i$ such that
\begin{equation}
b_k^j d a_j^i - a_j^i b_k^\ell \theta_\ell^j + \ti{\theta}_k^i = a^i_{k\ell} \ti{\theta}^\ell,
\end{equation}
which can be rewritten as
\begin{equation} \label{eqnda0}
da^i_m - a^i_j  \theta^j_m +  a^k_m \ti{\theta}^i_k= a^i_{k\ell} a^k_m \ti{\theta}^\ell.
\end{equation}
Note that $a^i_{k\ell} \ti{e}_i  \ti{\theta}^k  \ti{\theta}^\ell$ can be interpreted as the difference of the two connections $\displaystyle{\ti{\nabla} - \nabla}$.
Also, if $g$ and $\ti{g}$ are quasi-K\"ahler, from \eqref{eqnbda} we see that we have
$a^i_{k\ell} = a^i_{\ell k}$.
We will now calculate a formula for $\ti{\Delta} u$.  

\begin{lemma} \label{lemmaDeltau} For $g$ and $\ti{g}$ almost-Hermitian metrics, and $a^i_j$, $a^i_{k\ell}$, $b^i_j$ as defined above, we have
$$\frac{1}{2} \ti{\Delta} u =  a^i_{k\ell} \ov{a^i_{p\ell}} a^k_j \ov{a^p_j} - \ov{a^i_j} a^k_j \ti{R}^i_{k\ell\ov{\ell}} + \ov{a^i_j} a^i_r b^q_\ell \ov{b^s_\ell} R^r_{jq \ov{s}}.$$
\end{lemma}

\begin{proof1} 
Applying the exterior derivative to (\ref{eqnda0}), using the first and second structure equations and simplifying, we have
\begin{equation*}
\begin{split}
&- a^i_j \Psi^j_m + a^k_{j\ell} a^j_m \ti{\theta}^\ell \wedge \ti{\theta}^i_k  + a^k_m \ti{\Psi}^i_k  
=  a^k_m da^i_{k\ell} \wedge \ti{\theta}^\ell  \\
&- a^i_{k\ell} a^j_m \ti{\theta}^k_j \wedge \ti{\theta}^\ell + a^i_{k\ell} a^k_{jp} a^j_m \ti{\theta}^p \wedge \ti{\theta}^\ell - a^i_{k\ell} a^k_m \ti{\theta}^\ell_j \wedge \ti{\theta}^j + a^i_{k\ell} a^k_m \ti{\Theta}^\ell.
\end{split}
\end{equation*}
Multiplying by $b^m_r$ and rearranging, we obtain
\begin{eqnarray} \nonumber
\lefteqn{(da^i_{r\ell} +  a^i_{k\ell} a^k_{rj} \ti{\theta}^j+ a^k_{r\ell} \ti{\theta}^i_k   - a^i_{k\ell} \ti{\theta}^k_r - a^i_{rj} \ti{\theta}^j_\ell ) \wedge \ti{\theta}^\ell} \\ \label{eqndaw}
& & \qquad \qquad =  - b^m_r \Psi^j_m a^i_j  + \ti{\Psi}^i_r - a^i_{r\ell} \ti{\Theta}^\ell.
\end{eqnarray}
Define $a^i_{r\ell p}$ and $a^i_{r\ell\ov{p}}$ by
\begin{equation} \label{eqnda00}
da^i_{r\ell} +  a^i_{k\ell} a^k_{rj} \ti{\theta}^j+ a^k_{r\ell} \ti{\theta}^i_k   - a^i_{k\ell} \ti{\theta}^k_r - a^i_{rj} \ti{\theta}^j_\ell  = a^i_{r\ell p} \ti{\theta}^p + a^i_{r\ell \ov{p}} \ov{\ti{\theta}^p}.
\end{equation}
Then taking the (1,1) part  of (\ref{eqndaw}) we see that
\begin{equation} \label{eqnda000}
a^i_{r\ell \ov{p}} \ov{\ti{\theta}^p} \wedge \ti{\theta}^\ell = (- \ti{R}^i_{r\ell \ov{p}} + a_j^i b_r^m b_\ell^q \ov{b^s_p} R^j_{mq\ov{s}}) \ov{\ti{\theta}^p} \wedge \ti{\theta}^\ell,
\end{equation}
where we recall that by definition
\begin{eqnarray*}
 (\ti{\Psi}^i_r)^{(1,1)} & = & - \ti{R}^i_{r\ell \ov{p}} \ov{\ti{\theta}^p} \wedge \ti{\theta}^\ell \\
 (\Psi^j_m)^{(1,1)} & = & - R^j_{mq \ov{s}} \ov{\theta^s} \wedge \theta^q.
\end{eqnarray*}

Note that
\begin{eqnarray} \label{eqndu1}
du & = & \ov{a_j^i} da_j^i + a_j^i d \ov{a_j^i}.
\end{eqnarray}
Then we see that from (\ref{eqnda0}),
\begin{eqnarray} \nonumber
du & = & \ov{a_j^i} (a^i_{k\ell} a^k_j \ti{\theta}^\ell + a^i_m {\theta}^m_j - a_j^k \ti{\theta}_k^i ) + a_j^i ( \ov{a^i_{k\ell} a^k_j \ti{\theta}^\ell} + \ov{a^i_m\theta^m_j} - \ov{a^k_j\ti{\theta}^i_k}  ).\\ \label{eqndu0}
& = &  \ov{a_j^i} a^i_{k\ell} a^k_j \ti{\theta}^\ell + a_j^i \ov{a^i_{k\ell}a^k_j\ti{\theta}^\ell}.
\end{eqnarray}
Hence $\partial u = \ov{a_j^i} a^i_{k\ell} a^k_j \ti{\theta}^\ell$.    Applying the exterior derivative to this and 
substituting from (\ref{eqnda0}), (\ref{eqnda00}) and (\ref{eqnda000}) we have,
\begin{eqnarray*}
d \partial u & = & a^i_{k\ell} \ov{a^i_j} a^k_{pq} a^p_j \ti{\theta}^q \wedge \ti{\theta}^\ell + a^i_{k\ell} a^k_j \ov{a^i_{pq}a^p_j\ti{\theta}^q} \wedge \ti{\theta}^\ell \\ 
&& \mbox{} 
+ \ov{a^i_j} a^k_j (a^i_{k\ell p} \ti{\theta}^p + a^i_{k\ell\ov{p}} \ov{\ti{\theta}^p} - a^i_{r\ell} a^r_{kp} \ti{\theta}^p) \wedge \ti{\theta}^\ell + \ov{a^i_j} a^k_j a^i_{k\ell} \ti{\Theta}^\ell \\
& =  & a^i_{k\ell} a^k_j \ov{a^i_{pq}a^p_j\ti{\theta}^q} \wedge \ti{\theta}^\ell  + \ov{a^i_j} a^k_j a^i_{k\ell p} \ti{\theta}^p \wedge \ti{\theta}^\ell \\
&& \mbox{} + \ov{a^i_j} a^k_j (- \ti{R}^i_{k\ell \ov{p}} + a_r^i b_k^m b_\ell^q \ov{b^s_p} R^r_{mq\ov{s}}) \ov{\ti{\theta}^p} \wedge \ti{\theta}^\ell + \ov{a^i_j} a^k_j a^i_{k\ell} \ti{\Theta}^\ell.
\end{eqnarray*}
Hence
$$(d \partial u)^{(1,1)} = a^i_{k\ell} a^k_j \ov{a^i_{pq}a^p_j\ti{\theta}^q} \wedge \ti{\theta}^\ell - 
\ov{a^i_j} a^k_j \ti{R}^i_{k\ell \ov{p}} \ov{\ti{\theta}^p} \wedge \ti{\theta}^\ell + \ov{a^i_j} a_r^i  b_\ell^q \ov{b^s_p} R^r_{jq\ov{s}} \ov{\ti{\theta}^p} \wedge \ti{\theta}^\ell.$$
Then from the definition of the canonical Laplacian, we have proved the lemma.
\end{proof1}

Now let $\nu = \det (a_i^j)$ and set $v = |\nu|^2 = \nu \ov{\nu}$, which is the ratio of the volume forms of $\ti{g}$ and $g$. We have the following lemma.

\begin{lemma} \label{lemmaDeltalogv}
For $g$ and $\ti{g}$ almost-Hermitian metrics, and $v$ as above, the following identities hold.
\begin{enumerate}
\item[(i)] $\displaystyle{(d \partial \log v)^{(1,1)} = - R_{k \ov{l}} \theta^k \wedge \ov{\theta^l} + \ti{R}_{k \ov{l}} a^k_i \ov{a^l_j} \theta^i \wedge \ov{\theta^j}}$
\item[(ii)] $\displaystyle{\Delta \log v = 2R - 2  \ti{R}_{k \ov{l}} a_i^k \ov{a_i^l}.}$
%\item[(ii)] $\displaystyle{(d\partial v)^{(2,0)} = - v (K^i_{ikl}  - \ti{K}^i_{ipq} a^p_k a^q_l)  \theta^k \wedge \theta^l}$
\end{enumerate}
\end{lemma}
\begin{proof1}  This proof is essentially contained in  \cite{GH}, but we include it here for the reader's convenience.  
Write $\nu_j^i$ for the $(i,j)$th cofactor of the matrix $(a_i^j)$, so that $\nu_j^i=\nu b^i_j$.  Then
$$ d \nu = \nu_j^i da_i^j.$$
From (\ref{eqnda0}) we have
$$da^i_m - a^i_j \theta^j_m + a^k_m \ti{\theta}^i_k = a^i_{k\ell} a^k_m a^\ell_r \theta^r.$$
Hence
\begin{eqnarray} \nonumber
d\nu & = & \nu_j^i ( a^j_{pq} a^p_i a^q_k \theta^k+ a_k^j \theta_i^k - a_i^k \ti{\theta}^j_k) \\ \label{eqndnu1}
& = & \nu_k  \theta^k + \nu (\theta^i_i - \ti{\theta}^i_i),
\end{eqnarray}
for $\nu_k = \nu_j^i a^j_{pq} a^p_i a^q_k   $.  Now
\begin{eqnarray*}
dv & = & \ov{\nu} d\nu + \nu d \ov{\nu} \\
& = & \ov{\nu} ( \nu_k \theta^k + \nu ( \theta^i_i - \ti{\theta}^i_i)) + \nu ( \ov{\nu_k} \ov{\theta^k} + \ov{\nu} (\ov{\theta^i_i} - \ov{\ti{\theta}^i_i})) \\
& = & \ov{\nu} \nu_k \theta^k+ \nu \ov{\nu_k} \ov{\theta^k}.
\end{eqnarray*}
Hence 
$\partial v = \ov{\nu} \nu_k \theta^k.$  Define $v_k$ and $v_{\ov{k}}$ by $dv = v_k \theta^k + v_{\ov{k}} \theta^k$.  Then $v_k = \ov{\nu} \nu_k$.  
Applying the exterior derivative to (\ref{eqndnu1}) and using the second structure equation we have
\begin{eqnarray*}
0 & = & d( \nu_k \theta^k) + d \nu \wedge (\theta^i_i - \ti{\theta}^i_i) + \nu d(\theta^i_i - \ti{\theta}^i_i) \\
& = & d(\nu_k \theta^k) + \nu_k \theta^k \wedge (\theta^i_i - \ti{\theta}^i_i) + \nu (\Psi_i^i - \ti{\Psi}^i_i).
\end{eqnarray*}
Multiplying by $\ov{\nu}$ and using (\ref{eqndnu1}) again we have
\begin{eqnarray*}
0 & = & \ov{\nu} d (\nu_k \theta^k) + \nu_k \theta^k \wedge (\ov{\nu_\ell} \ov{\theta^\ell} - d \ov{\nu} ) + v (\Psi^i_i - \ti{\Psi^i_i}) \\
& = & d (\ov{\nu} \nu_k \theta^k) + \nu_k \ov{\nu_\ell} \theta^k \wedge \ov{\theta^\ell} + v (\Psi^i_i - \ti{\Psi}^i_i ).
\end{eqnarray*}
Consider the (1,1) part
\begin{eqnarray}\nonumber
(d \partial v)^{(1,1)} & = & - \nu_k \ov{\nu_\ell} \theta^k \wedge \ov{\theta^\ell} - v (\Psi^i_i - \ti{\Psi}^i_i)^{(1,1)} \\ \label{derivativev}
& = & - \frac{ v_k \ov{v_\ell}}{v} \theta^k \wedge \ov{\theta^\ell} - v R_{k \ov{\ell}} \theta^k \wedge \ov{\theta^\ell} + v \ti{R}_{k \ov{\ell}} a^k_i \ov{a^\ell_j} \theta^i \wedge \ov{\theta^j}.
\end{eqnarray}
We also have
$$d \partial \log v=\frac{d\partial v}{v}+\frac{\partial v\wedge\ov{\partial} v}{v^2},$$
which combines with \eqref{derivativev} to give (i).
From the definition of the canonical Laplacian we immediately obtain (ii).
\end{proof1}

We now return to the Calabi-Yau equation (\ref{eqnCY2}):
\begin{equation} \label{eqnCY1}
\det \ti{g} = e^{2F} \det g,
\end{equation}
for smooth $F$, where $g$ is almost-Hermitian and $\ti{g}$ is almost-K\"ahler.  Note that this equation can be rewritten in terms of $v$ as
\begin{equation} \label{eqnCY3}
\log v = F.
\end{equation}
We have the following lemma.

\begin{lemma} \label{lemmaDeltaulowerbound}  Suppose that $g$ is almost-Hermitian and $\ti{g}$ is quasi-K\"ahler and solves the Calabi-Yau equation (\ref{eqnCY2}).  Then
\begin{equation*}
\begin{split}
 \textit{(i)} & \quad\ti{\Delta} u = 2  a^i_{k\ell} \ov{a^i_{p\ell}} a^k_j \ov{a^p_j} + \Delta F - 2R  +8N^\ell_{\ov{p}\,\ov{i}}\ov{N^p_{\ov{\ell}\,\ov{i}}} + 2 \ov{a^p_i} a^p_j b^k_q \ov{b^\ell_q} \cur_{i\ov{j}k\ov{\ell}} \\
%\ov{a^i_j} a^i_r b^q_l \ov{b^s_l}R^r_{jq \ov{s}})\\
%& +  2\Delta F - 2R+8 \ov{a^i_j} a^k_j\left(\ti{N}^\ell_{\ov{p}\,\ov{i}}\ov{\ti{N}^p_{\ov{\ell}\,\ov{k}}}+\ti{N}^p_{\ov%{l}\,\ov{i}}\ov{\ti{N}^k_{\ov{p}\,\ov{l}}}\right)
 \textit{(ii)} & \quad \ti{\Delta} \log u   \geq \frac{1}{u} \Bigl( \Delta F - 2R  +8N^\ell_{\ov{p}\,\ov{i}}\ov{N^p_{\ov{\ell}\,\ov{i}}} + 2 \ov{a^p_i} a^p_j b^k_q \ov{b^\ell_q} \cur_{i\ov{j}k\ov{\ell}} \Bigr),   
\end{split}
\end{equation*}
where $\cur_{i\ov{j} k \ov{\ell}} = R^j_{ik\ov{\ell}} + 4N^r_{\ov{\ell} \, \ov{j}} \ov{N^i_{\ov{r} \, \ov{k}}}$.
\end{lemma}
\begin{proof1}
From Lemma \ref{lemmaDeltau}, Lemma \ref{lemmaDeltalogv} and the identity \eqref{riccicurv},
\begin{eqnarray*}
\quad\ti{\Delta} u & = & 2  (a^i_{k\ell} \ov{a^i_{p\ell}} a^k_j \ov{a^p_j} + \ov{a^i_j} a^i_r b^q_\ell \ov{b^s_\ell}R^r_{jq \ov{s}}) \\
&& \mbox{} +  \Delta F - 2R+8 \ov{a^i_j} a^k_j\left(\ti{N}^\ell_{\ov{p}\,\ov{i}}\ov{\ti{N}^p_{\ov{\ell}\,\ov{k}}}+\ti{N}^p_{\ov{\ell}\,\ov{i}}\ov{\ti{N}^k_{\ov{p}\,\ov{\ell}}}\right).
\end{eqnarray*}
Using \eqref{eqnN0}, we have
$$ \ov{a^i_j} a^k_j\left(\ti{N}^\ell_{\ov{p}\,\ov{i}}\ov{\ti{N}^p_{\ov{\ell}\,\ov{k}}}+\ti{N}^p_{\ov{\ell}\,\ov{i}}\ov{\ti{N}^k_{\ov{p}\,\ov{\ell}}}\right)= N^\ell_{\ov{p}\,\ov{i}}\ov{N^p_{\ov{\ell}\,\ov{i}}}+
\ov{a^k_s}a^k_j \ov{b^t_\ell} b^r_\ell N^p_{\ov{t}\,\ov{j}}\ov{N^s_{\ov{p}\,\ov{r}}},$$
giving (i).

For part (ii)
we compute,
\begin{eqnarray*}
\ti{\Delta} \log u & = & \frac{1}{u} \left( \ti{\Delta} u - \frac{| \ti{\nabla} u|_{\ti{g}}^2}{u} \right) \\
%& = & \frac{1}{u} \biggl( 2  (a^i_{kl} \ov{a^i_{pl}} a^k_j \ov{a^p_j} + \ov{a^i_j} a^i_r b^q_l \ov{b^s_l}R^r_%{jq \ov{s}})   +  2\Delta F - 2R -  \frac{| \ti{\nabla} u|_{\ti{g}}^2}{u}\\
%&& \mbox{} + 8N^l_{\ov{p}\,\ov{i}}\ov{N^p_{\ov{l}\,\ov{i}}}+
%8\ov{a^k_s}a^k_j \ov{b^t_l} b^r_l N^p_{\ov{t}\,\ov{j}}\ov{N^s_{\ov{p}\,\ov{r}}}\biggr) \\
& = & \frac{1}{u} \biggl( 2  a^i_{k\ell} \ov{a^i_{p\ell}} a^k_j \ov{a^p_j} + 8N^\ell_{\ov{p}\,\ov{i}}\ov{N^p_{\ov{\ell}\,\ov{i}}} + 2 \ov{a^p_i} a^p_j b^k_q \ov{b^\ell_q} \cur_{i\ov{j}k\ov{\ell}}   \\ 
&& \mbox{} +  \Delta F - 2R -  \frac{| \ti{\nabla} u|_{\ti{g}}^2}{u} \biggr).
\end{eqnarray*}
It remains to prove the inequality
\begin{equation} \label{eqnineq1}
| \ti{\nabla} u|_{\ti{g}}^2 \le  2u 
 a^i_{k\ell} \ov{a^i_{p\ell}} a^k_j \ov{a^p_j}.
\end{equation}
From (\ref{eqndu0}) we have 
$$ | \ti{\nabla} u |^2_{\ti{g}} = 2  u_\ell \ov{u_\ell},$$
where $u_\ell = \ov{a^i_j} a^i_{k\ell} a^k_j= \ov{a^i_j} B_{\ell j}^i$,  where $B^i_{\ell j} = a^i_{k\ell} a^k_j$.  Then using the Cauchy-Schwarz inequality,
\begin{eqnarray*}
 u_\ell \ov{u_\ell} & = & \sum_{i,j,\ell,p,q} \ov{a^i_j} B^i_{\ell j} a^p_q \ov{B^p_{\ell q}} \\
& \le  & \sum_{i,j,p,q} \left(\sum_\ell | \ov{a^i_j} B^i_{\ell j} |^2\right)^{1/2} \left( \sum_\ell |a^p_q \ov{B^p_{\ell q}}|^2 \right)^{1/2} \\
& = & \left( \sum_{i,j} \left( \sum_\ell |a^i_j|^2 |B^i_{\ell j}|^2 \right)^{1/2} \right)^2 \\
& = & \left( \sum_{i,j} |a^i_j| \left( \sum_\ell |B^i_{\ell j}|^2 \right)^{1/2} \right)^2 \\
& \le & \left( \sum_{i,j} |a^i_j|^2 \right) \left( \sum_{i,j,\ell} |B^i_{\ell j}|^2 \right) \\
& = & u  a^i_{k\ell} \ov{a^i_{p\ell}} a^k_j \ov{a^p_j},
\end{eqnarray*}
which gives (\ref{eqnineq1}).
\end{proof1}

Finally we can give the proof of Theorem \ref{theoremC2}.

\bigskip
\noindent
{\bf Proof of Theorem \ref{theoremC2}} \ 
Note that from the Calabi-Yau equation and the arithmetic-geometric means inequality, $u = \frac{1}{2} \tr{g}{\ti{g}}$ is bounded below away from zero by a positive constant depending only on $\sup|F|$.
Then from Lemma \ref{lemmaDeltaulowerbound} there exists $C'$ and $A'$ such that
$$\ti{\Delta} \log u \ge -C' - A' \tr{\ti{g}}{g},$$
with $A'$ depending only on the lower bound of $\cur_{i\ov{j} k \ov{l}}$, and $C'$ depending only on  $J$, $\sup |F|$, $\Delta F$ and $R$.
We apply the maximum principle to $(\log u - 2A' \varphi)$.  Suppose that the maximum of this quantity is achieved at a point $x_0$.  Then at this point, using (\ref{eqntiDeltaphi}),
\begin{eqnarray*}
0 \ge \ti{\Delta} (\log u - 2A' \varphi) \ge - C' + A' \tr{\ti{g}}{g} - 4A'n.
\end{eqnarray*}
Hence $$(\tr{\ti{g}}{g})(x_0) \le \frac{4A'n + C'}{A'}.$$
Using the inequality
$$\frac{\sum_{i=1}^n \lambda_i}{\prod_{i=1}^n\lambda_i}\leq \frac{1}{(n-1)!} \left(\sum_{i=1}^n \frac{1}{\lambda_i} \right)^{n-1},$$
%\leq\frac{(\sum_i\lambda_i)^{n-1}}
%{\prod_i\lambda_i},$$
that holds for any set of real numbers $\lambda_i>0$, and using the Calabi-Yau equation again,
we see that $u$ can be bounded from above in terms of $\tr{\ti{g}}{g}$ and so we obtain an estimate $$u(x_0) \le C''.$$
It follows that for any $x \in M$,
$$ \log u(x) - 2A' \varphi (x) \le \log C'' - 2A' \inf_M \varphi,$$
and the theorem is proved. 
\hfill$\square$\medskip

\begin{remark} \label{remarktheoremC2}
Notice that if we assume $\cur(g)>0$ in Theorem \ref{theoremC2}, then from Lemma \ref{lemmaDeltaulowerbound} we have
$$\ti{\Delta} \log u \ge -C' +A' u,$$
for some positive constant $A'$ and the maximum principle immediately gives $u\leq C.$
\end{remark}

\setcounter{equation}{0}
\setcounter{lemma}{0}
\setcounter{theorem}{0}
\addtocounter{section}{1}
\bigskip 
\bigskip \pagebreak[3]
\noindent
{\bf 4. First derivative estimate of $\ti{g}$}
\bigskip

In this section we give an estimate on the derivative of an almost-K\"ahler metric $\ti{g}$ 
 solving the Calabi-Yau equation (\ref{eqnCY2}).  This is a generalization of the third order estimate of \cite{Ya} (see also the recent preprint \cite{PSS} for a succinct proof of the parabolic version of this estimate).  Define 
$$S = \frac{1}{4} | \nabla \ti{g} |_{\ti{g}}^2,$$
where $\nabla$ is the canonical connection associated to $g$, $J$.   Then we have the following theorem.

\begin{theorem} \label{theoremS}
Let $\ti{g}$ be a solution of (\ref{eqnCY2}) and suppose that there exists a constant $K$ such that $$\sup_M (\emph{tr}_g \ti{g}) \le K.$$
Then there exists a constant $C_0$ depending only on $g$, $J$, $F$ and $K$ such that
$$S \le C_0.$$
\end{theorem}

Before we prove this theorem, we will need a number of lemmas.
\begin{lemma}
$S$ can be written as
\begin{equation} \label{eqnS}
S = a^i_{k\ell} \ov{a^i_{k\ell}},
\end{equation}
where $a^i_{k\ell}$ is defined by
\begin{equation} \label{eqndaagain}
da_m^i - a_j^i \theta^j_m + a_m^k \ti{\theta}^i_k = a^i_{k\ell} a^k_m \ti{\theta}^\ell.
\end{equation}
\end{lemma}
\begin{proof1}
To see (\ref{eqnS}) we calculate as follows:
\begin{eqnarray*}
\nabla \left( \ti{\theta}^i \otimes \overline{\ti{\theta}^i} \right) & = & \nabla (a^i_j \theta^j) \otimes \ov{\ti{\theta}^i} + \ti{\theta^i} \otimes \nabla ( \ov{a^i_j \theta^j}) \\
& = & (da^i_j) b^j_k \otimes \ti{\theta}^k \otimes \ov{\ti{\theta}^i} - a^i_j \theta^j_k b^k_\ell \otimes \ti{\theta}^\ell \otimes \ov{\ti{\theta}^i} \\
&& \mbox{} + \ov{(d a^i_j) b^j_k} \otimes \ti{\theta}^i \otimes \ov{\ti{\theta}^k} - \ov{a^i_j b^k_\ell \theta^j_k} \otimes \ti{\theta}^i \otimes \ov{\ti{\theta}^\ell} \\
& = & (da^\ell_j - a^\ell_r \theta^r_j) b^j_k \otimes \ti{\theta}^k \otimes \ov{\ti{\theta}^\ell} + \ov{(da^k_j - a^k_r \theta^r_j) b^j_\ell} \otimes \ti{\theta}^k \otimes \ov{\ti{\theta}^\ell} \\
& = & (a^\ell_{rs} a^r_j \ti{\theta}^s - a^r_j \ti{\theta}^\ell_r) b^j_k \otimes \ti{\theta}^k \otimes \ov{\ti{\theta}^\ell} \\
&& \mbox{} + \ov{( a^k_{rs} a^r_j\ti{\theta}^s - a^r_j \ti{\theta}^k_r) b_\ell^j}  \otimes \ti{\theta}^k \otimes \ov{\ti{\theta}^\ell} \\
& = & a^\ell_{ks} \ti{\theta}^s \otimes \ti{\theta}^k \otimes \ov{\ti{\theta}^\ell} + \ov{a^k_{\ell s}} \, \ov{\ti{\theta}^s} \otimes \ti{\theta}^k \otimes \ov{\ti{\theta}^\ell}.
\end{eqnarray*}
Then since $\ti{g} = \ti{\theta}^i \otimes \ov{\ti{\theta}^i} + \ov{\ti{\theta}^i} \otimes \ti{\theta}^i$,  (\ref{eqnS}) follows immediately.
\end{proof1}

The following lemma gives a general formula for the Laplacian of $S$.

\begin{lemma} \label{lemmalapS} We have
\begin{eqnarray} \nonumber
\frac{1}{2} \ti{\Delta} S
 & = & \left| a^i_{k\ell p} - a^i_{r\ell} a^r_{kp} \right|_{\ti{g}}^2 + |a^i_{k\ell\ov{p}}|_{\ti{g}}^2 + 
\ov{a^i_{k\ell}} a^i_{r\ell} \ti{R}^r_{kp\ov{p}} + \ov{a^i_{k\ell}} a^i_{kj} \ti{R}^j_{\ell p\ov{p}} - \ov{a^i_{k\ell}} a^r_{k\ell} \ti{R}^i_{r p \ov{p}}
  \\ \nonumber
&& \mbox{} +2 \mathrm{Re} \biggl( \ov{a^i_{k\ell}} \biggl( b^m_k b^q_\ell \ov{b^s_p} R^j_{mq\ov{s}} a^i_{r p} a^r_j 
- a^i_j b^q_\ell \ov{b^s_p} R^j_{mq\ov{s}} a^r_{kp} b^m_r   \\ \nonumber
&&  \mbox{}  - a^i_j b^m_k \ov{b^s_p} R^j_{mq\ov{s}} a^r_{\ell p} b^q_r
 + a^i_j b^m_k b^q_\ell \ov{b^s_p} b_p^u  R^j_{mq\ov{s},u} - \ti{R}_{k \ov{i},\ell} \\\nonumber && 
+4 \ti{N}^p_{\ov{q}\,\ov{i},\ell}\ov{\ti{N}^q_{\ov{p}\,\ov{k}}}+4 \ti{N}^p_{\ov{q}\,\ov{i}}\ov{\ti{N}^q_{\ov{p}\,\ov{k},\ov{\ell}}}
+4\ti{N}^p_{\ov{q}\,\ov{i},\ell}\ov{\ti{N}^k_{\ov{p}\,\ov{q}}}
+ 4\ti{N}^p_{\ov{q}\,\ov{i}}\ov{\ti{N}^k_{\ov{p}\,\ov{q},\ov{\ell}}}\\ 
&&+4\ti{N}^i_{\ov{q}\,\ov{p},k}\ov{\ti{N}^q_{\ov{p}\,\ov{\ell}}}
 + 2 \ov{\ti{N}^k_{\ov{\ell} \, \ov{p}, ip}} \biggr)\biggr). \label{eqnlemmalapS}
 \end{eqnarray}
\end{lemma}
\begin{proof1}
First, recall from (\ref{eqndaw}) and (\ref{eqnda00}) that $a^i_{r\ell p}$ and $a^i_{r\ell\ov{p}}$ are defined by
\begin{equation} \label{eqnda00again}
da^i_{r\ell} +  a^i_{k\ell} a^k_{rj} \ti{\theta}^j+ a^k_{r\ell} \ti{\theta}^i_k   - a^i_{k\ell} \ti{\theta}^k_r - a^i_{rj} \ti{\theta}^j_\ell  = a^i_{r\ell p} \ti{\theta}^p + a^i_{r\ell \ov{p}} \ov{\ti{\theta}^p},
\end{equation}
and that
\begin{eqnarray}  
(a^i_{r\ell p} \ti{\theta}^p + a^i_{r\ell \ov{p}} \ov{\ti{\theta}^p}) \wedge \ti{\theta}^\ell  \label{eqndawagain}
&= & - b^m_r \Psi^j_m a^i_j  + \ti{\Psi}^i_r - a^i_{r\ell} \ti{\Theta}^\ell.
\end{eqnarray}
Define functions $a^i_{r\ell p,q}, a^i_{r\ell p,\ov{q}}, a^i_{r\ell\ov{p},q},$ and $a^i_{r\ell\ov{p},\ov{q}}$ by the formulas
\begin{equation}\label{deriva1}
da^i_{r\ell p}+a^k_{r\ell p}\ti{\theta}^i_k -a^i_{r\ell q}\ti{\theta}^q_p    -a^i_{k\ell p}\ti{\theta}^k_r-a^i_{rjp}\ti{\theta}^j_\ell=a^i_{r\ell p,q}\ti{\theta}^q+
a^i_{r\ell p,\ov{q}}\ov{\ti{\theta}^q},
\end{equation}
\begin{equation}\label{deriva2}
da^i_{r\ell\ov{p}}+a^k_{r\ell\ov{p}}\ti{\theta}^i_k -a^i_{r\ell\ov{q}}\ov{\ti{\theta}^q_p}   -
a^i_{k\ell\ov{p}}\ti{\theta}^k_r-a^i_{rj\ov{p}}\ti{\theta}^j_\ell=a^i_{r\ell\ov{p},q}\ti{\theta}^q+
a^i_{r\ell\ov{p},\ov{q}}\ov{\ti{\theta}^q}.
\end{equation}
Applying the exterior derivative to (\ref{eqnda00again}), using the last two definitions, and canceling many terms we get
\begin{eqnarray} \nonumber
\lefteqn{a^i_{r\ell p,q} \ti{\theta}^q \wedge \ti{\theta}^p + a^i_{r\ell p,\ov{q}} \ov{\ti{\theta}^q}\wedge \ti{\theta}^p +a^i_{rl\ov{p},q} \ti{\theta} \wedge \ov{\ti{\theta}^p} + a^i_{rl\ov{p}, \ov{q}} \ov{\ti{\theta}^q} \wedge \ov{\ti{\theta}^p}} \\ \nonumber && \mbox{}    + 
a^i_{r\ell p} \ti{\Theta}^p + a^i_{r\ell\ov{p}} \ov{\ti{\Theta}^p}  \\ \nonumber
& = & \mbox{} - a^k_{rp} a^i_{s\ell} a^s_{kt} \ti{\theta}^t \wedge \ti{\theta}^p + a^k_{rp} a^i_{k\ell s} \ti{\theta}^s \wedge \ti{\theta}^p + a^k_{rp} a^i_{k\ell\ov{s}} \ov{\ti{\theta}^s} \wedge \ti{\theta}^p \\ \nonumber
& &\mbox{}  - a^i_{k\ell} a^k_{sp} a^s_{rt} \ti{\theta}^t \wedge \ti{\theta}^p + a^i_{k\ell} a^k_{rpt} \ti{\theta}^t \wedge \ti{\theta}^p + a^i_{k\ell} a^k_{rp\ov{t}} \ov{\ti{\theta}^t} \wedge \ti{\theta}^p \\
&& \mbox{} +a^i_{k\ell} a^k_{rp} \ti{\Theta}^p - a^i_{k\ell} \ti{\Psi}^k_r - a^i_{rp} \ti{\Psi}^p_\ell + a^k_{r\ell} \ti{\Psi}^i_k, \label{eqnda0000}
\end{eqnarray}
which will be useful later.  To calculate the canonical Laplacian of $S$ with respect to $\ti{g}$,  first compute 
\begin{eqnarray*}
\partial S & = & \ov{a^i_{k\ell}} \partial a^i_{k\ell} + a^i_{k\ell} \ov{ \ov{\partial} a^i_{k\ell}} \\
& = & \bigl( \ov{a^i_{k\ell}} a^i_{k\ell p}  + a^i_{k\ell} \ov{a^i_{k\ell\ov{p}}} - \ov{a^i_{k\ell}} a^i_{r\ell} a^r_{kp} \bigr) \ti{\theta}^p.
\end{eqnarray*}
Then compute
\begin{eqnarray}\nonumber
d (\partial S) &=& \bigl(a^i_{k\ell p}\ov{a^i_{k\ell q}\ti{\theta}^q}+a^i_{k\ell p}\ov{a^i_{k\ell\ov{q}}}\ti{\theta}^q-a^i_{k\ell p}\ov{a^i_{r\ell}a^r_{kq}\ti{\theta}^q}
+\ov{a^i_{k\ell}}a^i_{k\ell p,q}\ti{\theta}^q\\ \nonumber
&&\mbox{} +\ov{a^i_{k\ell}}a^i_{k\ell p,\ov{q}}\ov{\ti{\theta}^q}+\ov{a^i_{k\ell\ov{p}}}a^i_{k\ell q}\ti{\theta}^q
+\ov{a^i_{k\ell\ov{p}}}a^i_{k\ell\ov{q}}\ov{\ti{\theta}^q}-\ov{a^i_{k\ell\ov{p}}}a^i_{r\ell}a^r_{kq}\ti{\theta}^q\\ \nonumber
&&\mbox{} +a^i_{k\ell}\ov{a^i_{k\ell\ov{p},q}\ti{\theta}^q}+a^i_{k\ell}\ov{a^i_{k\ell\ov{p},\ov{q}}}\ti{\theta}^q
-\ov{a^i_{k\ell q}}a^i_{r\ell}a^r_{kp}\ov{\ti{\theta}^q}-\ov{a^i_{k\ell\ov{q}}}a^i_{r\ell}a^r_{kp}\ti{\theta}^q\\ \nonumber
&&\mbox{} +a^i_{r\ell}a^r_{kp}\ov{a^i_{j\ell}a^j_{kq}\ti{\theta}^q}-\ov{a^i_{k\ell}}a^r_{kp}a^i_{r\ell q}\ti{\theta}^q
-\ov{a^i_{k\ell}}a^r_{kp}a^i_{r\ell\ov{q}}\ov{\ti{\theta}^q}\\ 
&&\mbox{} +\ov{a^i_{k\ell}}a^r_{kp}a^i_{j\ell}a^j_{rq}\ti{\theta}^q-\ov{a^i_{k\ell}}a^i_{r\ell}a^r_{kpq}\ti{\theta}^q
-\ov{a^i_{k\ell}}a^i_{r\ell}a^r_{kp\ov{q}}\ov{\ti{\theta}^q}\\ \nonumber
&&\mbox{} +\ov{a^i_{k\ell}}a^i_{r\ell}a^r_{jp}a^j_{kq}\ti{\theta}^q\bigr)\wedge\ti{\theta}^p+\bigl( \ov{a^i_{k\ell}} a^i_{k\ell p}  + a^i_{k\ell} \ov{a^i_{k\ell\ov{p}}} - \ov{a^i_{k\ell}} a^i_{r\ell} a^r_{kp}\bigr) \ti{\Theta}^p,
\end{eqnarray}
and hence
\begin{eqnarray} \label{eqnddS1}\nonumber
\left( d( \partial S) \right)^{(1,1)} &=&
\bigl(a^i_{k\ell p}\ov{a^i_{k\ell q}}-a^i_{k\ell p}\ov{a^i_{r\ell}a^r_{kq}}+\ov{a^i_{k\ell}}a^i_{k\ell p,\ov{q}}
+\ov{a^i_{k\ell\ov{p}}}a^i_{k\ell\ov{q}}\\ \nonumber
&&\mbox{} +a^i_{k\ell}\ov{a^i_{k\ell\ov{p},q}}
-\ov{a^i_{k\ell q}}a^i_{r\ell}a^r_{kp} +a^i_{r\ell}a^r_{kp}\ov{a^i_{j\ell}a^j_{kq}}\\ 
&&\mbox{} 
-\ov{a^i_{k\ell}}a^r_{kp}a^i_{r\ell\ov{q}}
-\ov{a^i_{k\ell}}a^i_{r\ell}a^r_{kp\ov{q}}\bigr)\ov{\ti{\theta}^q}\wedge\ti{\theta}^p.
\end{eqnarray}

Then taking the (1,1) part of (\ref{eqnda0000}) we see that
\begin{eqnarray}\label{eqndbara}\nonumber
a^i_{k\ell p,\ov{q}}\ov{\ti{\theta}^q}\wedge\ti{\theta}^p&=&\bigl(a^i_{k\ell\ov{q},p}+a^i_{r\ell\ov{q}}a^r_{kp}+a^i_{r\ell}a^r_{kp\ov{q}}\\
&&\mbox{}
+a^i_{r\ell}\ti{R}^r_{kp\ov{q}}+a^i_{kj}\ti{R}^j_{\ell p\ov{q}} -a^r_{k\ell}\ti{R}^i_{rp\ov{q}}  \bigr)\ov{\ti{\theta}^q}\wedge\ti{\theta}^p.
\end{eqnarray}
Multiplying \eqref{eqndbara} by $\ov{a^i_{k\ell}}$, substituting into \eqref{eqnddS1} and using the formula for the Laplacian, we obtain
\begin{eqnarray*}
\frac{1}{2} \ti{\Delta}  S & = & a^i_{k\ell p} \ov{a^i_{k\ell p}} - a^i_{k\ell p} \ov{a^i_{r\ell} a^r_{kp}} + \ov{a^i_{k\ell\ov{p}}} a^i_{k\ell\ov{p}} - a^i_{r\ell} a^r_{kp} \ov{a^i_{k\ell p}}  \\
&& \mbox{} + a^i_{r\ell} a^r_{kp} \ov{a^i_{s\ell} a^s_{kp}} - \ov{a^i_{k\ell}} a^r_{kp} a^i_{r\ell\ov{p}} - \ov{a^i_{k\ell}} a^i_{r\ell} a^r_{kp\ov{p}} \\
&& \mbox{} + \ov{a^i_{k\ell}} a^i_{k\ell \ov{p},p} + \ov{a^i_{k\ell}} a^r_{kp} a^i_{r\ell\ov{p}} + \ov{a^i_{k\ell}} a^i_{r\ell} a^r_{kp\ov{p}} \\
&& \mbox{} + \ov{a^i_{k\ell}} a^i_{r\ell} \ti{R}^r_{kp\ov{p}} + \ov{a^i_{k\ell}} a^i_{kj} \ti{R}^j_{\ell p\ov{p}} - \ov{a^i_{k\ell}} a^r_{k\ell} \ti{R}^i_{r p \ov{p}} + a^i_{k\ell} \ov{a^i_{k\ell\ov{p},p}} \\
& = & a^i_{k\ell p} \ov{a^i_{k\ell p}} + \ov{a^i_{k\ell\ov{p}}} a^i_{k\ell\ov{p}} - 2 \textrm{Re} (a^i_{k\ell p} \ov{a^i_{r\ell} a^r_{kp}})  \\
&& \mbox{} + a^i_{r\ell} a^r_{kp} \ov{a^i_{s\ell} a^s_{kp}} + \ov{a^i_{k\ell}} a^i_{r\ell} \ti{R}^r_{kp\ov{p}} + \ov{a^i_{k\ell}} a^i_{kj} \ti{R}^j_{\ell p\ov{p}} - \ov{a^i_{k\ell}} a^r_{k\ell} \ti{R}^i_{r p \ov{p}} \\
&& \mbox{} +2 \textrm{Re} ( \ov{a^i_{k\ell}} a^i_{k\ell \ov{p},p}).
\end{eqnarray*}
Completing the square, we obtain
\begin{eqnarray} \nonumber
\frac{1}{2} \ti{\Delta} S & = & \left| a^i_{k\ell p} - a^i_{r\ell} a^r_{kp} \right|_{\ti{g}}^2 + |a^i_{k\ell\ov{p}}|_{\ti{g}}^2 + 
\ov{a^i_{k\ell}} a^i_{r\ell} \ti{R}^r_{kp\ov{p}} + \ov{a^i_{k\ell}} a^i_{kj} \ti{R}^j_{\ell p\ov{p}}
 \\ \label{eqnlapS1}
&& \mbox{}  - \ov{a^i_{k\ell}} a^r_{k\ell} \ti{R}^i_{r p \ov{p}} +2 \textrm{Re} ( \ov{a^i_{k\ell}} a^i_{k\ell \ov{p},p}).
\end{eqnarray}
To calculate the last term, take the (1,1) part of (\ref{eqndawagain}) to obtain
\begin{equation}\label{deriva3}
a^i_{k\ell \ov{p}} = a^i_j b^m_k b^q_\ell \ov{b^s_p} R^j_{mq\ov{s}} - \ti{R}^i_{k\ell\ov{p}}.
\end{equation}
Now recall from (\ref{eqnda0}) that
\begin{equation} \label{eqnderiva}
da^i_m - a^i_j \theta^j_m + a^k_m \ti{\theta}^i_k = a^i_{k\ell} a^k_m \ti{\theta}^\ell.
\end{equation}
Similarly we have
\begin{equation}\label{derivb}
db^j_k+b^r_k\theta^j_r-b^j_i\ti{\theta}^i_k=-b^j_i a^i_{k\ell}\ti{\theta}^\ell.
\end{equation}
Taking the exterior derivative of \eqref{deriva3}, using \eqref{deriva1}, \eqref{deriva2}, \eqref{eqnderiva} and \eqref{derivb} we get
\begin{eqnarray} \nonumber
a^i_{k\ell\ov{p},t}\ti{\theta}^t+a^i_{k\ell\ov{p},\ov{t}}\ov{\ti{\theta}^t}&=&a^i_j b^m_k b^q_\ell \ov{b^s_p} R^j_{mq\ov{s},u}\theta^u
+a^i_j b^m_k b^q_\ell \ov{b^s_p} R^j_{mq\ov{s},\ov{u}}\ov{\theta^u}\\ \nonumber
&&\mbox{} -\ti{R}^i_{k\ell\ov{p},t}\ti{\theta}^t-\ti{R}^i_{k\ell\ov{p},\ov{t}}\ov{\ti{\theta}^t}
+b^m_k b^q_\ell \ov{b^s_p}a^r_j a^i_{rt}R^j_{mq\ov{s}}\ti{\theta}^t\\ \nonumber
&&\mbox{} -b^m_r b^q_\ell \ov{b^s_p}a^i_j a^r_{kt}R^j_{mq\ov{s}}\ti{\theta}^t
-b^m_k b^q_r \ov{b^s_p} a^i_j a^r_{\ell t}R^j_{mq\ov{s}}\ti{\theta}^t\\ 
&&\mbox{} 
-b^m_k b^q_\ell \ov{b^s_r} a^i_j \ov{a^r_{pt}}R^j_{mq\ov{s}}\ov{\ti{\theta}^t},
\end{eqnarray}
whose $(1,0)$ part gives
\begin{eqnarray} \nonumber
a^i_{k\ell\ov{p},t} & = & b^m_k b^q_\ell \ov{b^s_p} R^j_{mq\ov{s}} a^i_{r t} a^r_j - a^i_j b^q_\ell \ov{b^s_p} R^j_{mq\ov{s}} a^r_{kt} b^m_r - a^i_j b^m_k \ov{b^s_p} R^j_{mq\ov{s}} a^r_{\ell t} b^q_r \\ \label{eqncalcda}
&& \mbox{} + a^i_j b^m_k b^q_\ell \ov{b^s_p}  b^u_t R^j_{mq\ov{s},u} -  \ti{R}^i_{k\ell\ov{p},t}.
\end{eqnarray}
Now from \eqref{commutation3},\eqref{commutation2} and \eqref{riccicurv}
\begin{eqnarray} \nonumber
\ti{R}^i_{k\ell\ov{p},p} & = & \ti{R}^i_{kp\ov{p},\ell} - 2 \ti{K}^i_{kp\ell, \ov{p}}-4\ti{K}^i_{k\ov{q}\,\ov{p}}\ov{\ti{N}^q_{\ov{p}\,\ov{\ell}}}\\
 \nonumber  \label{eqncalcdR}
& = & \ti{R}_{k \ov{i},\ell}-4 \ti{N}^p_{\ov{q}\,\ov{i},\ell}\ov{\ti{N}^q_{\ov{p}\,\ov{k}}}-4 \ti{N}^p_{\ov{q}\,\ov{i}}\ov{\ti{N}^q_{\ov{p}\,\ov{k},\ov{\ell}}}
-4\ti{N}^p_{\ov{q}\,\ov{i},\ell}\ov{\ti{N}^k_{\ov{p}\,\ov{q}}} \\
&&\mbox{}- 4\ti{N}^p_{\ov{q}\,\ov{i}}\ov{\ti{N}^k_{\ov{p}\,\ov{q},\ov{\ell}}}- 2 \ti{K}^i_{kp\ell, \ov{p}}-4\ti{N}^i_{\ov{q}\,\ov{p},k}\ov{\ti{N}^q_{\ov{p}\,\ov{\ell}}},
\end{eqnarray}
and using \eqref{commutation3} again we see that
\begin{eqnarray}
\ti{K}^i_{kp\ell, \ov{p}}  & = & \ov{\ti{N}^k_{\ov{\ell} \, \ov{p}, ip}} \label{eqncalcaK}.
\end{eqnarray}
Combining (\ref{eqnlapS1}), (\ref{eqncalcda}), (\ref{eqncalcdR}) and (\ref{eqncalcaK}) gives 
(\ref{eqnlemmalapS}).
\end{proof1}

To deal with the terms involving derivatives of $\ti{N}^i_{\ov{j} \, \ov{k}}$ in (\ref{eqnlemmalapS}) we need another lemma.

\pagebreak[3]
\begin{lemma} \label{lemmaN} We have
\begin{enumerate}
\item[(i)] $\displaystyle{\ti{N}^i_{\ov{j} \, \ov{k}, m} = \ov{b^r_j b^s_k} b^\ell_m a^i_t N^t_{\ov{r} \, \ov{s}, \ell} + \ov{b^r_j b^s_k} a^\ell_t N^t_{\ov{r} \, \ov{s}} a^i_{\ell m}}$
\item[(ii)] $\displaystyle{\ti{N}^i_{\ov{j} \, \ov{k}, \ov{m}} = \ov{b^r_j b^s_k b^\ell_m} a^i_t N^t_{\ov{r} \, \ov{s}, \ov{\ell}} -
 \ov{b^r_\ell b^s_k} a^i_t N^t_{\ov{r} \, \ov{s}} \ov{a^\ell_{jm}}-\ov{b^r_j b^s_\ell} a^i_t N^t_{\ov{r} \, \ov{s}} \ov{a^\ell_{km}}}$
\item[(iii)] $\displaystyle{ \left| a^i_{k\ell}\ti{N}^k_{\ov{\ell} \, \ov{p}, ip} \right|_{\ti{g}} \le C (S + 1) + \frac{1}{2} \left| a^i_{k\ell p} - a^i_{r\ell} a^r_{kp}\right|^2_{\ti{g}}}$,
\end{enumerate}
for a constant $C$ depending only on $g$, $J$, $\sup_M \emph{tr}_g {\ti{g}}$ and $\sup_M \emph{tr}_{\ti{g}} g$.
\end{lemma}
\begin{proof1}
Recall from \eqref{eqnN0} that we have
$$\ti{N}^i_{\ov{j} \, \ov{k}} =  \ov{b^r_j b^s_k} a^i_t N^t_{\ov{r} \, \ov{s}}.$$
Applying the exterior derivative to this and using \eqref{eqnderiva}, \eqref{derivb} and \eqref{derivativeN} we obtain
\begin{eqnarray}\nonumber
\lefteqn{\ti{N}^i_{\ov{j}\, \ov{k}, m}\ti{\theta}^m+\ti{N}^i_{\ov{j}\, \ov{k}, \ov{m}}\ov{\ti{\theta}^m} } \\ \nonumber & =&
\ov{b^r_j b^s_k}b^p_m a^i_t N^{t}_{\ov{r}\, \ov{s}, p}\ti{\theta}^m+\ov{b^r_j b^s_k b^p_m} a^i_t N^{t}_{\ov{r}\, \ov{s}, \ov{p}}\ov{\ti{\theta}^m}+
a^i_{\ell m}  \ov{b^r_j b^s_k} a^{\ell}_t N^t_{\ov{r} \, \ov{s}}  \ti{\theta}^m\\ \label{eqniandii}
&&\mbox{} -  \ov{a^{\ell}_{jm} b^r_{\ell} b^s_k} a^i_t  N^t_{\ov{r} \, \ov{s}}   \ov{\ti{\theta}^m}
-  \ov{a^{\ell}_{km} b^r_j b^s_{\ell}} a^i_t N^t_{\ov{r} \, \ov{s}}  \ov{\ti{\theta}^m}.
\end{eqnarray}
Equating the $(1,0)$ and $(0,1)$ parts of \eqref{eqniandii} gives (i) and (ii).
For (iii), apply the exterior derivative to (i) and substitute from (\ref{eqnda00again}) to get
\begin{eqnarray}\label{calcolo}\nonumber
\ti{N}^i_{\ov{j}\,\ov{k},mp}&=&\ov{b^r_j b^s_k}b^\ell_m b^q_p a^i_t N^t_{\ov{r}\,\ov{s},\ell q}+\ov{b^r_j b^s_k}b^q_p a^\ell_t a^i_{\ell m} N^t_{\ov{r}\,\ov{s},q}
-\ov{b^r_j b^s_k}b^\ell_q a^i_t a^q_{mp}N^t_{\ov{r}\,\ov{s},\ell}\\
&&\mbox{}+\ov{b^r_j b^s_k}b^\ell_m a^q_t a^i_{qp}N^t_{\ov{r}\,\ov{s},\ell}
+\ov{b^r_j b^s_k}a^\ell_t N^t_{\ov{r}\,\ov{s}} a^i_{\ell mp}.
\end{eqnarray}
The only term that is not comparable to $\sqrt{S}$ is the last one.
To deal with this we first compute, using \eqref{commutation1}, \eqref{switchr} and \eqref{bianchin}
\begin{equation*}
\begin{split}
\ti{R}^i_{jk\ov{\ell}}&=\ti{R}^\ell_{kj\ov{i}}+4 \ti{N}^i_{\ov{p}\,\ov{\ell}}\ov{\ti{N}^p_{\ov{j}\,\ov{k}}}+
4\ti{N}^p_{\ov{i}\,\ov{\ell}}\ov{\ti{N}^k_{\ov{p}\,\ov{j}}}\\
&= \ti{R}^\ell_{jk\ov{i}}+ 4 \ti{N}^\ell_{\ov{p}\,\ov{i}}\ov{\ti{N}^p_{\ov{k}\,\ov{j}}}+4 \ti{N}^i_{\ov{p}\,\ov{\ell}}\ov{\ti{N}^p_{\ov{j}\,\ov{k}}}+4\ti{N}^p_{\ov{i}\,\ov{\ell}}\ov{\ti{N}^k_{\ov{p}\,\ov{j}}}\\
&=\ti{R}^\ell_{jk\ov{i}}+ 4 \ti{N}_{\ov{\ell}\,\ov{p}\,\ov{i}}\ov{\ti{N}_{\ov{p}\,\ov{k}\,\ov{j}}}+4 \ti{N}_{\ov{i}\,\ov{\ell}\,\ov{p}}\ov{\ti{N}_{\ov{p}\,\ov{k}\,\ov{j}}}+4\ti{N}_{\ov{p}\,\ov{i}\,\ov{\ell}}\ov{\ti{N}_{\ov{p}\,\ov{k}\,\ov{j}}}
+4\ti{N}_{\ov{p}\,\ov{i}\,\ov{\ell}}\ov{\ti{N}_{\ov{j}\,\ov{p}\,\ov{k}}}\\
&=\ti{R}^\ell_{jk\ov{i}}+4\ti{N}^p_{\ov{i}\,\ov{\ell}}\ov{\ti{N}^j_{\ov{p}\,\ov{k}}},
\end{split}
\end{equation*}
and use this, \eqref{commutation3}, \eqref{commutation2} and \eqref{curvident} to compute
\begin{eqnarray}\label{calcolo2}\nonumber
2\ov{\ti{N}^k_{\ov{\ell} \, \ov{p}, ip}}&=&2\ti{K}^i_{kp\ell,\ov{p}} \\ \nonumber
& = & 4\ov{\ti{K}^k_{iqp}\ti{N}^q_{\ov{p}\,\ov{\ell}}}+
\ti{R}^i_{kp\ov{p},\ell}-\ti{R}^i_{k\ell\ov{p},p} \\ \nonumber
& =& 4\ov{\ti{K}^k_{iqp}\ti{N}^q_{\ov{p}\,\ov{\ell}}}+
\ti{R}^p_{kp\ov{i},\ell}
-\ti{R}^p_{k\ell\,\ov{i},p}
+4\ti{N}^q_{\ov{i}\,\ov{p},\ell}\ov{\ti{N}^k_{\ov{q}\,\ov{p}}}+4\ti{N}^q_{\ov{i}\,\ov{p}}\ov{\ti{N}^k_{\ov{q}\,\ov{p},\ov{\ell}}} \\ \nonumber && \mbox{}
-4\ti{N}^q_{\ov{i}\,\ov{p},p}\ov{\ti{N}^k_{\ov{q}\,\ov{\ell}}}-4\ti{N}^q_{\ov{i}\,\ov{p}}\ov{\ti{N}^k_{\ov{q}\,\ov{\ell},\ov{p}}}\\ \nonumber
&=&2\ov{\ti{N}^k_{\ov{\ell} \, \ov{p}, pi}}+4\ti{N}^i_{\ov{p} \, \ov{q}, k} \ov{\ti{N}^q_{\ov{p}\,\ov{\ell}}}-4\ti{N}^p_{\ov{i}\, \ov{q}, k} \ov{\ti{N}^q_{\ov{p}\,\ov{\ell}}}
+4\ti{N}^q_{\ov{i}\,\ov{p},\ell}\ov{\ti{N}^k_{\ov{q}\,\ov{p}}}+4\ti{N}^q_{\ov{i}\,\ov{p}}\ov{\ti{N}^k_{\ov{q}\,\ov{p},\ov{\ell}}}\\
&&\mbox{}-4\ti{N}^q_{\ov{i}\,\ov{p},p}\ov{\ti{N}^k_{\ov{q}\,\ov{\ell}}}-4\ti{N}^q_{\ov{i}\,\ov{p}}\ov{\ti{N}^k_{\ov{q}\,\ov{\ell},\ov{p}}}.
\end{eqnarray}
This means that, up to an error comparable to $\sqrt{S}$, we can interchange the last two covariant
derivatives on $\ti{N}$. Finally recall from \eqref{eqnbda} that
$$a^i_{k\ell}\ti{\theta}^\ell\wedge\ti{\theta}^k=a^i_j T^j_{pq}b^p_k b^q_\ell\ti{\theta}^\ell\wedge\ti{\theta}^k,$$
and so 
\begin{equation}\label{scambio}
a^i_{k\ell}=a^i_{\ell k}+2a^i_j b^p_k b^q_\ell T^j_{pq}.
\end{equation}
From \eqref{calcolo2}, \eqref{calcolo}, \eqref{eqnN0} and \eqref{scambio},
\begin{equation} \nonumber
\begin{split}
\left| a^i_{k\ell}\ti{N}^k_{\ov{\ell} \, \ov{p}, ip} \right|_{\ti{g}}& \leq C(S+1)+\left| a^i_{k\ell}\ti{N}^k_{\ov{\ell} \, \ov{p}, pi} \right|_{\ti{g}}\\
&\leq C(S+1)+\left|a^i_{k\ell}\ti{N}^q_{\ov{\ell}\,\ov{p}} a^k_{qpi}\right|_{\ti{g}}\\
&\leq C(S+1)+\left|a^i_{k\ell}\ti{N}^q_{\ov{\ell}\,\ov{p}} (a^k_{qpi}-a^k_{rp}a^r_{qi})\right|_{\ti{g}}+\left|a^i_{k\ell}a^k_{rp}a^r_{qi}\ti{N}^q_{\ov{\ell}\,\ov{p}}\right|_{\ti{g}}\\
&\leq C(S+1)+\frac{1}{2}\left|a^k_{qpi}-a^k_{rp}a^r_{qi}\right|^2_{\ti{g}}+\left|a^i_{k\ell}a^k_{rp}a^r_{qi}\ti{N}^q_{\ov{\ell}\,\ov{p}}\right|_{\ti{g}}\\
&\leq C(S+1)+\frac{1}{2}\left|a^k_{qpi}-a^k_{rp}a^r_{qi}\right|^2_{\ti{g}}+\left|a^i_{k\ell}a^k_{rp}a^r_{iq}\ti{N}^q_{\ov{\ell}\,\ov{p}}\right|_{\ti{g}},
\end{split}
\end{equation}
where the constant $C$ differs from line to line, and 
where we have used the inequality 
$$2ab \le {\varepsilon}a^2 + \frac{1}{\varepsilon}b^2,$$
for any $\epsilon>0$ and any real numbers $a$ and $b$.
Finally, using \eqref{bianchin} we can see that the term $a^i_{k\ell}a^k_{rp}a^r_{iq}\ti{N}_{\ov{q}\,\ov{\ell}\,\ov{p}}$ vanishes:
\begin{equation*}
\begin{split}
a^i_{k\ell}a^k_{rp}a^r_{iq}\ti{N}_{\ov{q}\,\ov{\ell}\,\ov{p}}&=\frac{1}{3}(a^i_{k\ell}a^k_{rp}a^r_{iq}+a^i_{kp}a^k_{rq}a^r_{i\ell}+a^i_{kq}a^k_{r\ell}a^r_{ip})
\ti{N}_{\ov{q}\,\ov{\ell}\,\ov{p}}\\
&=\frac{1}{3}a^i_{k\ell}a^k_{rp}a^r_{iq}(\ti{N}_{\ov{q}\,\ov{\ell}\,\ov{p}}+\ti{N}_{\ov{p}\,\ov{q}\,\ov{\ell}}+\ti{N}_{\ov{\ell}\,\ov{p}\,\ov{q}})=0.
\end{split}
\end{equation*}
\end{proof1}

We can now prove the following lemma.

\begin{lemma} \label{lemmalowerboundlapS}
Let $\ti{g}$ be an almost-K\"ahler metric solving the Calabi-Yau equation (\ref{eqnCY2}) and suppose that there exists a constant $K$ such that $$\sup_M (\emph{tr}_g \ti{g}) \le K.$$  Then there exist constants $C_1$, $C_2$ depending only on $g$, $J$, $F$ and $K$ such that
\begin{equation} \label{eqnineqlapS}
\ti{\Delta} S \ge -C_1 S -C_2.
\end{equation}
\end{lemma}
\begin{proof1}
By assumption, the $a^i_j$ and $b^i_j$ are uniformly bounded.  From \eqref{eqnCY3} and \eqref{eqndpartialf} we have
$$(d \partial \log v)^{(1,1)} = -F_{p \ov{q}} \theta^p \wedge \ov{\theta^q}.$$
Then from Lemma \ref{lemmaDeltalogv}, we have
$$\ti{R}_{k\ov{\ell}} = - F_{p \ov{q}} b^p_k \ov{b^q_\ell}  + R_{p\ov{q}} b^p_k \ov{b^q_\ell}.$$
It follows that $|\ti{R}_{k\ov{\ell}}|_{\ti{g}}^2 \le C$ and $|\ti{R}_{k\ov{\ell},p}|_{\ti{g}}^2 \le C(S+1)$, for a constant $C$ depending only on $g$, $J$, $F$ and $K$.  Then the inequality (\ref{eqnineqlapS}) follows from Lemma \ref{lemmalapS} and Lemma \ref{lemmaN}.
\end{proof1}

Finally, we complete the proof of Theorem \ref{theoremS}.

\bigskip
\noindent
{\bf Proof of Theorem \ref{theoremS}} \ Following \cite{Ya} we apply the maximum principle to $S + C'u$, for a constant $C'$ to be determined later.  Note that from Lemma \ref{lemmaDeltaulowerbound} (i), we have
$$\ti{\Delta} u \ge C_3^{-1} S - C_4,$$
for positive constants $C_3$ and $C_4$ depending only on $g$, $J$, $F$ and $K$.  Choose $C' = C_3 (C_1+1)$ then from Lemma \ref{lemmalowerboundlapS} we see that
$$\ti{\Delta} (S+C'u) \ge S - C_2 - C' C_4,$$
and then by the maximum principle $S$ is bounded from above by $C_0=C_2 + C'C_4+ C'K$.
\hfill$\square$\medskip

\setcounter{equation}{0}
\setcounter{lemma}{0}
\addtocounter{section}{1}
\bigskip
\bigskip \pagebreak[3]
\noindent
{\bf 5. Proof of Theorem 1}
\bigskip

Let $\ti{g}$ solve the Calabi-Yau equation (\ref{eqnCY2}).  We will write $\nabla_g$ and $dV_g$ for the Levi-Civita covariant derivative and volume form associated to the metric $g$.  We have the following lemma.

\begin{lemma}  \label{lemmaalpha} For every $\alpha>0$ there exists a constant $C$ depending only on $(M, \Omega, J)$, $F$ and $\alpha$ such that
$$- \inf_M \varphi \le C + \log \left( \int_M e^{-\alpha \varphi} dV_g \right)^{1/\alpha}.$$
\end{lemma}

\bigskip
\noindent
{\bf Proof of Lemma \ref{lemmaalpha}} \ 
Let $\delta>0$ be a small constant.  In the following $C$ will denote a uniform constant, depending only on $\delta$ and the fixed data, which may change from line to line.
Define $w=e^{-B\varphi}$ for $B = \frac{1}{1-\delta}A$.  Write $\gamma=1-\delta>0$.
For $p\ge 1$, from Theorem \ref{theoremC2} and the Calabi-Yau equation,
\begin{eqnarray*}
\int_M | \nabla_g w^{p/2}|^2 dV_g & \le & - C \int_M (\tr{{g}}{\ti{g}}) d e^{- \frac{Bp\varphi}{2}} \wedge J d e^{- \frac{Bp\varphi}{2}} \wedge \ti{\omega}^{n-1} \\
& \le & -C p^2 e^{-B\gamma \inf_M \varphi}  \int_M e^{-B(p-\gamma) \varphi} d \varphi \wedge J d \varphi \wedge \ti{\omega}^{n-1} \\
& = &  C\frac{p^2}{p-\gamma}  \| w \|_{C^0}^{\ga}  \int_M d \left( e^{-B(p-\gamma) \varphi} \right) \wedge J d \varphi \wedge \ti{\omega}^{n-1} \\
& \le & C p \| w \|_{C^0}^{\ga} \int_M   w^{p-\gamma}  \ti{\Delta} \varphi \, \ti{\omega}^n\\
& \le & Cp \| w \|_{C^0}^{\gamma}  \int_M   w^{p-\gamma}  dV_g,
\end{eqnarray*}
using the fact that $\ti{\Delta} \varphi \le 2n$ from \eqref{eqntiDeltaphi}.  The Sobolev inequality gives us, for $\beta =\frac{n}{n-1}$,
$$\| f \|_{2 \beta}^2 \le C ( \|\nabla_g f \|_2^2 + \|f \|_2^2),$$
where $\| \ \|_q$ denotes the $L^q$ norm with respect to $g$ (we allow later $0<q<1$, defined in the obvious way).  Applying this to $f=w^{p/2}$, we obtain
\begin{eqnarray*}
\left( \int_M w^{p \beta} dV_g \right)^{1/\beta} & \le & C \left( \int_M | \nabla_g w^{p/2}|^2 dV_g + \int_M w^p dV_g \right) \\
& \le & C p \| w \|_{C^0}^{\gamma} \int_M w^{p-\gamma} dV_g.
\end{eqnarray*}
Raising to the power $1/p$ we have
$$\| w \|_{p\beta} \le C^{1/p} p^{1/p} \| w \|_{C^0}^{\gamma/p} \| w \|_{p-\ga}^{(p-\ga)/p}.$$
By the same iteration as in \cite{We1} we replace $p$ with $p\beta + \gamma$ to obtain for $k=1,2, \ldots$,
$$ \| w \|_{p_k \beta} \le C(k) \| w \|_{C^0}^{1-a(k)} \| w \|_{p-\gamma}^{a(k)},$$
where
\begin{eqnarray*}
p_k & = & p \beta^k + \gamma (1+ \beta + \beta^2 + \cdots \beta^{k-1}) \\
C(k) & = & C^{(1+ \beta + \cdots + \beta^k)/p_k}p_0^{\beta^k/p_k} p_1^{\beta^{k-1}/p_k} \cdots p_k^{1/p_k} \\
a(k) & = & \frac{(p-\gamma)\beta^k}{p_k}.
\end{eqnarray*}
Set $p=1$ and let $k \rightarrow \infty$.  Since $C(k)$ is uniformly bounded from above and $a(k) \rightarrow a \in (0,1)$, we have
$$\| w \|_{C^0} \le C \| w \|_{\delta},$$
and choosing $\delta$ sufficiently small completes the proof of the lemma.
\hfill$\square$\medskip

From this lemma, Theorem \ref{theoremC2} and Theorem \ref{theoremS} we have the estimate
$$\|\ti{g}\|_{C^1}\leq C,$$
where $C$ depends on $\Omega, J,\sigma, \alpha$ and $I_\alpha(\varphi)$. It remains to prove the higher order estimates.  Following \cite{We2}, define
a $1$-form $a$ by the equations
$$\ti{\omega}=\Omega-\frac{1}{2} d(Jd\varphi)+da,$$
and $d^*_{\ti{g}}a=0$, where $d^*_{\ti{g}}$ is the formal adjoint of $d$ associated to $\ti{g}$. Note that $a$ is defined
only up to the addition of a harmonic $1$-form. From the definition of $\varphi$ it follows that $da\wedge\ti{\omega}^{n-1}=0$.
Let's call $\mathcal{P}:\Lambda^2(M)\to\Lambda^2(M)$ the map that associates to a $2$-form $\gamma$ its $(2,0)+(0,2)$ part, so that
$$\mathcal{P}\gamma(X,Y)=\frac{1}{2}(\gamma(X,Y)-\gamma(JX,JY)).$$
Since $\ti{\omega}$ is compatible with $J$ we have $\mathcal{P}\ti{\omega}=0$, but in general $\mathcal{P}\Omega\neq 0$.
Now set $f=\varphi$ in \eqref{eqnddu} and take the $(2,0)$ part to get
$$\varphi_{ij}\theta^j\wedge\theta^i+\varphi_{\ov{k}}\ov{N^k_{\ov{j}\,\ov{i}}}\theta^j\wedge\theta^i + \varphi_k T^k_{ji} \theta^j \wedge \theta^i =0.$$
Applying $\mathcal{P}$ to \eqref{derivf},
\begin{eqnarray*}
\mathcal{P} d(Jd\varphi)&=&2\sqrt{-1}\mathcal{P}d\partial\varphi \\
& = & 2\sqrt{-1}\left(\varphi_{ij}\theta^j\wedge\theta^i+ \varphi_k T^k_{ji} \theta^j \wedge \theta^i +
\varphi_k N^k_{\ov{j}\,\ov{i}}\ov{\theta^j}\wedge\ov{\theta^i}\right)\\
&=&2\sqrt{-1} \left( \varphi_k N^k_{\ov{j}\,\ov{i}} \ov{\theta^j} \wedge \ov{\theta^i} -\varphi_{\ov{k}}\ov{N^k_{\ov{j}\,\ov{i}}} \theta^j \wedge \theta^i\right),  
\end{eqnarray*}
which  involves only one derivative of $\varphi$. Now the $1$-form $a$ satisfies the following system
\begin{equation}\label{system}
\left\{  \begin{array}{rcl}
 da\wedge\ti{\omega}^{n-1}&=&0 \\
\mathcal{P}da &=&-\mathcal{P}\Omega+ \sqrt{-1} \left( \varphi_k N^k_{\ov{j} \, \ov{i}} \ov{\theta^j} \wedge \ov{\theta^i} - \varphi_{\ov{k}} \ov{N^k_{\ov{j} \, \ov{i}}} \theta^j \wedge \theta^i \right)\\
d^*_{\ti{g}}a&=&0,
\end{array}
\right.
\end{equation}
which is elliptic (its symbol is injective, although not invertible if $n>2$).   

Note that the kernel of \eqref{system} consists of the harmonic 1-forms.  Indeed, $da \wedge \ti{\omega}^{n-1}=0$ and $\mathcal{P}(da)=0$ together imply that $*da = - c_n \ti{\omega}^{n-2} \wedge da$ for some universal constant $c_n$.  Then if $a$ is in the kernel of \eqref{system}, we have  $\| da \|^2_{L^2(\ti{g})}=0$ after integrating by parts.  Since $d^*_{\ti{g}}a=0$, we see that $a$ is harmonic with respect to $\ti{g}$.

Fix any $0<\beta<1$. Since $\ti{g}$ is uniformly bounded in $C^\beta$, we can apply the elliptic Schauder estimates to 
\eqref{eqntiDeltaphi} to get a bound $\|\varphi\|_{C^{2+\beta}}\leq C$. Hence the right hand side of \eqref{system} is bounded in $C^{1+\beta}$, and
the coefficients of the system have a $C^\beta$ bound, so assuming that $a$ is orthogonal to the harmonic 1-forms, the elliptic estimates applied to \eqref{system} give $C^{2+\beta}$ bounds on $a$.
By differentiating the Calabi-Yau equation in a direction $\partial/\partial x^i$ we obtain
\begin{equation}\label{equationderiv}
\ti{\Delta}(\partial_i\varphi)+\{\textrm{lower order terms}\}=2\partial_i F+g^{pq}\partial_i g_{pq},
\end{equation}
where the lower order terms may contain up to two derivatives of $\varphi$ or $a$, and so are bounded in $C^\beta$.
Applying the Schauder estimates again we get $\|\varphi\|_{C^{3+\beta}}\leq C$, and using \eqref{system} again we get $\| a\|_{C^{3+\beta}}\leq C$.
Now a bootstrapping argument using \eqref{equationderiv} and \eqref{system} gives the required higher order estimates. This completes the proof of Theorem 1.
\hfill$\square$\medskip

\setcounter{equation}{0}
\setcounter{lemma}{0}
\addtocounter{section}{1}
\bigskip
\bigskip \pagebreak[3]
\noindent
{\bf 6. Proof of Theorem 2}
\bigskip

As before, let $\ti{g}$ be an almost-K\"ahler metric solving (\ref{eqnCY2}).  Let $g$ be an almost-Hermitian metric 
with the property that $\cur(g) \ge 0$.  By the argument of the last section, to prove Theorem 2, it suffices to prove a uniform upper bound for $u= \frac{1}{2} \tr{g}{\ti{g}}$.

From Lemma \ref{lemmaDeltaulowerbound}, we have
$$\ti{\Delta} u \ge - C,$$
for a constant $C$ depending only on the fixed data.  We claim that this is enough to bound $u$ uniformly from above.  Indeed, for $p\ge 1$,
\begin{eqnarray*}
\int_M | \nabla_g u^{p/2}|^2 dV_g  & \le & - C \int_M u du^{p/2} \wedge J d u^{p/2} \wedge \ti{\omega}^{n-1}\\
& = & - Cp^2 \int_M u^{p-1} du \wedge Jdu \wedge \ti{\omega}^{n-1} \\
& = & - Cp \int_M d (u^p) \wedge J du \wedge \ti{\omega}^{n-1} \\
& = & Cp \int_M u^p d (J du) \wedge \ti{\omega}^{n-1} \\
& = & - C p \int_M u^p (\ti{\Delta} u) \ti{\omega}^n \\
& \le & C p \int_M u^p dV_g.
\end{eqnarray*}
Hence
$$\int_M | \nabla u^{p/2}|^2 dV_g \le C p \int_M u^p dV_g.$$
Then from the Sobolev inequality, we obtain
$$\| u \|_{L^{p\beta}} \le C^{1/p} p^{1/p} \| u \|_{L^p},$$
for $\beta=\frac{n}{n-1}$.
Replacing $p$ with $p\beta$, iterating, and then setting $p=1$ we obtain
$$\| u \|_{C^0} \le C \| u \|_{L^1}.$$
But this last quantity is bounded, because from \eqref{cohomolo} and the Calabi-Yau equation \eqref{eqnCYK},
$$\int_M u dV_g\leq C\int_M \frac{\ti{\omega}^{n-1}\wedge\Omega}{\ti{\omega}^n}\Omega^n\leq C\int_M \ti{\omega}^{n-1}\wedge\Omega=C[\Omega]^n.$$
This completes the proof of Theorem 2.
\hfill$\square$\medskip

\pagebreak[3]
\bigskip
\noindent
{\bf Acknowledgements.} \ The first author would like to thank Yanir Rubinstein, Aleksandar Suboti\'c and 
Chen-Yu Chi for some helpful discussions.  The second author thanks Simon Donaldson for some useful and encouraging conversations, and his former advisor, D.H. Phong for his support and advice.  The second author is supported in part by NSF grant DMS 0504285.

\end{document}